\documentclass[a4paper,11pt,draft]{amsart}  
\usepackage{amssymb}   
\usepackage{amsmath,amsthm,geometry,graphicx} 
\usepackage{amscd}
\usepackage{psfrag,xspace} 
\usepackage[all,dvips,cmtip]{xy} 

\normalfont  
  
\newcommand{\RedefinitSymbole}[1]{%
  \expandafter\let\csname old\string#1\endcsname=#1 \let#1=\relax 
\newcommand{#1}{\csname old\string#1\endcsname\,}%
} 
\RedefinitSymbole{\forall} \RedefinitSymbole{\exists} 

\def\cf{\textit{cf.}\kern.3em} 
 
\def\ie{\textit{i.e.}\ } 
\def\resp{\textit{resp.}\kern.3em} 
\renewcommand{\k}{\kern2pt} 
 
\numberwithin{equation}{section} 
 
\DeclareMathOperator{\lcm}{lcm}

\DeclareMathOperator{\age}{age} 
  
\DeclareMathOperator{\ev}{ev} 
 \DeclareMathOperator{\orb}{orb}

\DeclareMathOperator{\diag}{diag} 
 \DeclareMathOperator{\pgcd}{gcd}

\DeclareMathOperator{\id}{id}

 \DeclareMathOperator{\codim}{codim} 
  
\renewcommand{\ker}{\mathop{\rm 
    Ker}\nolimits} \newcommand{\reg}{\mathrm{reg}} 
 
\newtheorem{prop}[equation]{Proposition} 
\newtheorem{thm}[equation]{Theorem}
\newtheorem{cor}[equation]{Corollary} 
\newtheorem{lem}[equation]{Lemma} 
\theoremstyle{definition} 
\newtheorem{defi}[equation]{Definition}
 
\newtheorem{rem}[equation]{\textbf{Remark}} 
\newtheorem{notation}[equation]{\textbf{Notation}} 
\newtheorem{expl}[equation]{\textbf{Example}}

\newcommand{\NN}{\mathbb{N}} 
\newcommand{\ZZ}{\mathbb{Z}} 
\newcommand{\QQ}{\mathbb{Q}} 
\newcommand{\CC}{\mathbb{C}} 
\newcommand{\PP}{\mathbb{P}} 
\newcommand{\ds}{\displaystyle} 
 
\newcommand{\bs}{\boldsymbol} 
\newcommand{\gr}{\mathrm{gr}}

\begin{document}
 \title{Orbifold quantum cohomology of
  weighted projective spaces}
 \author{Etienne Mann}
 \address{Etienne
  Mann\\ SISSA\\via Beirut 2-4\\ I-34014 TRIESTE\\ ITALY}
\email{mann@sissa.it} \urladdr{http://www-irma.u-strasbg.fr/$\tilde{\ 
  }$mann/}
\date{\today}

\begin{abstract}
In this article, we prove the following results.
\begin{itemize}
  \item We show a mirror theorem : the Frobenius manifold associated
    to the orbifold quantum cohomology of weighted projective space is
    isomorphic to the one attached to a specific Laurent polynomial,
  \item We show a reconstruction theorem, that is, we can
  reconstruct in an algorithmic way the full genus $0$ Gromov-Witten
  potential from the $3$-point invariants.
  \end{itemize}
\end{abstract}
\maketitle
 
 
 
\section{Introduction}\label{sec:Introduction}

Motivated by the works of  physicists E.\kern2pt Witten, R.\k
Dijkgraaf, E.\kern2pt Verlinde and H.\kern2pt Verlinde, B.\kern2pt
Dubrovin defined in   \cite{Dtft}  a Frobenius  structure on a complex
manifold.
Frobenius manifolds are complex manifolds endowed with a flat metric
and a product on the tangent bundle which
satisfies some compatibility conditions.

 In $2001$, S.\k Barannikov showed in  \cite{Bms} that the  Frobenius manifold coming
 from the quantum cohomology of the complex projective
 space of dimension $n$ is isomorphic to the 
  Frobenius manifold associated to the Laurent polynomial $x_{1}+
 \ldots+x_{n}+1/x_{1}\ldots x_{n}$.

 The goal of this article is to generalize this correspondence to
 weighted projective spaces. For this purpose we use the theory of orbifolds
 and the related constructions. In \cite{CRogw} and \cite{CRnco}, W.\k
 Chen and Y.\k Ruan define the orbifold cohomology ring via the
 orbifold Gromov-Witten invariants.  The orbifold cup product is
 defined as the degree zero part of the orbifold quantum product
 and one computes it via the Euler class of an obstruction bundle.
 The orbifold quantum product is defined by the Gromov-Witten
 potential.  So, as for manifolds, the orbifold quantum cohomology is
 naturally endowed with a Frobenius structure.
 
 On the other side, A.\k Douai and C.\k Sabbah (\cf \cite{DSgm1})
 explained how to build a canonical Frobenius manifold on the base
 space of a universal unfolding for any Laurent polynomial which is
 convenient and non-degenerate with respect to its Newton polyhedron.
 In particular, in \cite{DSgm2}, the authors described explicitly this
 construction for the polynomial $w_{0}u_{0}+ \cdots+ w_{n}u_{n}$
 restricted to $U:=\{(u_{0}, \ldots ,u_{n})\in \CC^{n+1}\mid \prod_{i}
 u_{i}^{w_{i}}=1\}$ where $w_{0}, \ldots ,w_{n}$ are positive integers
 which are relatively prime.
 
 In this article, we compare the Frobenius structures, whose existence
 is provided by the general results recalled above, on the orbifold
 quantum cohomology of the weighted projective space $\PP(w_{0},
 \ldots ,w_{n})$ (A side) and the one attached to the Laurent
 polynomial $f(u_{0}, \ldots ,u_{n}):=u_{0}+\cdots + u_{n}$ restricted
 to $U$ (B side).

 First we show a correspondence between "classical limits". To state
 this result, we need to introduce some notations. For the A side, we
 denote by $H^{2\star}_{\orb}(\PP(w_{0}, \ldots ,w_{n}),\CC)$ the
 orbifold cohomology of $\PP(w_{0}, \ldots ,w_{n})$, $\cup$ the
 orbifold cup product and $\langle\cdot,\cdot\rangle$ the orbifold
 Poincar\'e duality. For the B side, we consider the vector space
 $\Omega^{n}(U)/df\wedge\Omega^{n-1}(U)$ where $\Omega^{n}(U)$ is the
 space of algebraic $n$-forms on $U$. It is naturally endowed with an
 increasing filtration, called the Newton filtration and denoted by
 $\mathcal{N}_{\bullet}$, and a non-degenerate bilinear form. The choice of a
 volume form on $U$ gives us an identification of this vector space
 with the Jacobian ring of $f$. Hence, we get a product on this vector
 space.  As the product and the non-degenerate bilinear form respect
 the filtration $\mathcal{N}_{\bullet}$, we have a product, denoted by $\cup$,
 and a non-degenerate bilinear form, denoted by
 $[\![g]\!](\cdot,\cdot)$, on the graded space of
 $\Omega^{n}(U)/df\wedge\Omega^{n-1}(U)$ with respect to the Newton
 filtration.  The following theorem is shown in Section
 \ref{subsec:proof-class-corr}.
   
\begin{thm}[Classical correspondence] 
\label{thm:correspondance,classique} We have an isomorphism   
of graded Frobenius algebras between
$\left(H^{2\star}_{\orb}(\PP(w),\CC),\cup,\langle\cdot,\cdot\rangle\right)$
and
$\left(\gr^{\mathcal{N}}_{\star}\left(\Omega^{n}(U)/df\wedge\Omega^{n-1}(U)\right),\cup,[\![g]\!](\cdot,\cdot)\right).$
\end{thm} 

Note that, in a more general and algebraic context, A.\k Borisov, L.\k
Chen and G.\k Smith \cite{BCSocdms} computed the orbifold cohomology
ring for toric Deligne-Mumford stacks. We will not use these results
because, firstly we will use the techniques developed by W.\k Chen and
Y.\k Ruan and, secondly the author did not find in the literature a
complete and explicit description of weighted projective spaces as
toric Deligne-Mumford stacks.

Afterward, using \cite{CCLTsqcwps}, we prove two
propositions\footnote{In a previous version of this article, these propositions were conjectures.} (\cf
\ref{prop:invariant,nul} and \ref{prop:invariants,durs}) on the value
of some orbifold Gromov-Witten invariants with $3$ marked points and
we show in Section \ref{subsec:proof-quant-corr} that these
propositions imply an isomorphism between the Frobenius manifolds
coming from the A side and from the B side.
Let us note that Theorem \ref{thm:cond,init} shows that we can
reconstruct, in an algorithmic way, the full genus $0$ Gromov-Witten
invariants from the $3$-point invariants. This result is similar to
the first reconstruction theorem of M.\k Kontsevich and Y.\k Manin in
\cite[Theorem $3.1$]{KMgwqceg}.

The article is organized as follows. The first section is devoted to
Frobenius manifolds.  In the second section, we compute the orbifold
cohomology ring of weighted projective spaces. In the third section,
we compute the value of some specific Gromov-Witten invariants.  In
the fourth part, we briefly recall the results about the Laurent
polynomial $f:U\to \CC$.  In the last section, we give the proofs of
the two correspondences : the "classical correspondence" and the
isomorphism between the two Frobenius manifolds. 

\textbf{Acknowledgments :} I want to thank Claude Sabbah who
gives me such a nice subject for my thesis. His advices
were always relevant. I am also grateful to Claus Hertling who follows
my work during these years. I am indebted to Barbara Fantechi for
helpful discussions on algebraic stacks and for her interests in my work.

\section{Recalls on Frobenius manifolds}\label{sec:short-revi-frob}

Let $M$ be a complex manifold endowed with
\begin{itemize}
\item a perfect pairing $g: TM\times TM \to \CC$,
\item an associative and commutative product $\star$ on the complex
  tangent bundle $TM$ with unit $e$,
\item a vector field $\mathfrak{E}$, called the  Euler vector field. 
\end{itemize}
These data $(M,\star,e,g,\mathfrak{E})$ defined a Frobenius structure
on $M$ if they satisfy some compatibility conditions. We will not
write them because we will not use them explicitly. The reader can
find these conditions in Lecture $1$ of \cite{Dtft} (see also \cite[p.19]{Mfm},
\cite[p.146]{Hfm}, \cite[p.240]{Sdivf}).  Assume that $M$ is
simply-connected. Let $(t_{0}, \ldots ,t_{n})$ be a system of flat
coordinates on $M$.  According to Lemma $1.2$ in Lecture $1$ of
\cite{Dtft} (see also Section VII.$2.b$ in \cite{Sdivf}), there exists
a holomorphic function, called potential, $F:M\rightarrow \CC$ such
that for any $i,j,k$ in $\{1, \ldots ,n\}$, we have
    \begin{align*}  
    \frac{\partial^{3}F}{\partial t_{i}\partial t_{j}\partial t_{k}}=g(\partial_{t_{i}}\star\partial_{t_{j}},\partial_{t_{k}}).   
    \end{align*}  
    The potential is determined up to a polynomial of degree $2$.  As
    the product $\star$ is associative, the potential satisfies the
    WDVV equations. We have the following theorem.
 
\begin{thm}[\cite{Dtft}, lecture $3$ 
  ; see also \cite{Sdivf} p.$250$ or more generally Theorem $4.5$ in \cite{HMumc}] \label{thm:iso,frobenius,dubrovin}
  Let $g^{\circ}:\CC^{\mu}\times\CC^{ \mu}\rightarrow \CC$ be a
  perfect pairing.  Let $A_{0}^{\circ}$ be a semi-simple and regular
  matrix of size $\mu\times\mu$ such that
  $(A_{0}^{\circ})^{\ast}=A_{0}^{\circ}$. Let $A_{\infty}$ be a matrix
  of size $\mu\times\mu$ such that $A_{\infty}+A_{\infty}^{
    \ast}=k \cdot\id$ with $k\in \ZZ$.  Let $e^{\circ}$ be an
  eigenvector of $A_{\infty}$ for the eigenvalue $q$ such that
  $(e^{\circ}, A_{0}^{\circ}e^{\circ}, \ldots ,
  {A_{0}^{\circ}}^{\mu-1}e^{\circ})$ is a basis of $\CC^{\mu}$.  The
  data $(A^{\circ}_{0},A_{\infty},g^{\circ},e^{\circ})$ determined a
  unique germ of Frobenius manifold $((M,0),\star,e,g,\mathfrak{E})$ 
  such that via the isomorphism between $T_{0}M$ and $\CC^{\mu}$ we
  have   $g^{\circ}={g}(0)$, $A^{\circ}_{0}=\mathfrak{E}\star$,
  $A_{\infty}=(q+1)\id-\nabla{\mathfrak{E}}$ and $e^{\circ}=e(0)$.
\end{thm} 

In order to show an isomorphism between the Frobenius manifold coming
form $\PP(w)$ and the one associated to the Laurent polynomial $f$, we
will show that their initial conditions satisfy the hypothesis of the
theorem above and that they are equal.

\subsection{The A side} 
\label{sec:cote,A} 
 
We construct the Frobenius manifold on the complex vector space
$H^\star_{\orb}(\PP(w_{0}, \ldots ,w_{n}),\CC)$ of dimension  
$\mu:=w_{0}+\cdots+w_{n}$. The perfect pairing is the orbifold
Poincar\'e duality, denoted by 
$\langle \cdot,\cdot \rangle$.  In Section  
\ref{subsec:Orbif-cohom-weight}, we will define a basis $(\eta_{0}, \ldots ,\eta_{\mu-1})$
of the vector space $H^\star_{\orb}(\PP(w_{0}, \ldots ,w_{n}),\CC)$.
Denote by  
$(t_{0},\ldots,t_{\mu-1})$ the coordinates  
$H^\star_{\orb}(\PP(w_{0}, \ldots ,w_{n}),\CC)$ in this basis.  
The Euler field is defined by the following formula  
\begin{align*} 
  \mathfrak{E}&:=\mu\partial_{t_{1}}+\sum_{i=0}^{\mu-1}(1-\deg(\eta_{i})/2)t_{i}\partial_{t_{i}}. 
\end{align*} 
The big quantum product, denoted by $\star$, is defined with the full
Gromov-Witten potential of genus $0$, denoted by  $F^{GW}$, by the
following formula 
\begin{align*} 
 \frac{\partial^3 F^{GW}(t_{0},\ldots,t_{\mu-1})}{\partial 
 t_{i}\partial t_{j} \partial t_{k}}&=\langle \partial t_{i}\star 
 \partial t_{j}, \partial t_{k}\rangle  
\end{align*} 
The initial conditions of the Frobenius manifold are the data 
$(A_{0}^{\circ},A_{\infty},\langle\cdot,\cdot\rangle,\eta_{0})$ where 
$A^\circ_{0}:=\mathfrak{E}\star\mid_{\mathbf{t}=0}$ and 
$A_{\infty}:=\id-\nabla \mathfrak{E}$.

The matrix $A_{\infty}$ is easy to compute (see Proposition
\ref{prop:matrice,Ainfty}), but in order to compute the matrix
$A^\circ_{0}$, we have to compute the orbifold cup product (\cf
Section \ref{subsec:Orbif-cohom-ring}) and some specific Gromov-Witten
invariants with $3$ marked points (\cf Section
\ref{subsec:calcul-de-certains}). Via the correspondence, Theorem
\ref{thm:cond,init} implies that  we
can reconstruct the big quantum cohomology from the small one. In
particular, the proof of Theorem \ref{thm:cond,init} gives an
algorithm to do so.

\subsection{The B side}
 
Let  $U:=\{(u_{0}, \ldots ,u_{n})\in \CC^{n+1}\mid \prod_{i} 
u_{i}^{w_{i}}=1\}$. 
In the article \cite{DSgm2}, the  polynomial is $w_{0}u_{0}+\cdots+w_{n}u_{n}$ 
restricted to $U$ and the weights are relatively prime. 
In our case, we consider general weights and the polynomial 
$f$ is $u_{0}+\cdots+u_{n}$ restricted to $U$. 
Nevertheless, we will use the same techniques to show the following theorem.  
 
\begin{thm}[see Theorem \ref{thm:cond,init,sing}] 
 There exists a canonical Frobenius structure on any germ
  of universal unfolding of the Laurent polynomial $f(u_{0}, \ldots
  ,u_{n})=u_{0}+\cdots+u_{n}$ restricted to $U$.
\end{thm}

 
\numberwithin{equation}{section} 
  \section{Orbifold cohomology ring of weighted projective spaces}
\label{sec:Orbif-cohom-ring-}
In this section, we will describe explicitly the orbifold cohomology
ring of weighted projective spaces.

In the first part, we define the orbifold structure that we will
consider on weighted projective spaces.
In the second part, we give a natural $\CC$-basis of the orbifold
cohomology then we compute the orbifold Poincar\'e duality in this basis.
In the last part, we compute the orbifold cup product and we express
it in the basis defined in the second part. The obstruction
bundle is computed in Theorem \ref{thm:fibre,obstruction}.

In this article, we will use the following notations.  Let $n$ and
$w_{0},\ldots,w_{n}$ be some integers greater or equal to one. 

\subsection{Orbifold structure on weighted projective spaces}
\label{subsec:Orbif-struct-}
In this part, we describe the weighted projective spaces as
Deligne-Mumford stacks in Section \ref{subsubsec:weight-proj-spac-stacks} and
as orbifold in Section \ref{subsubsec:weight-proj-spac}.

\subsubsection{Weighted projective spaces as Deligne-Mumford stacks}
\label{subsubsec:weight-proj-spac-stacks}
We define the action of the multiplicative group $\CC^{\star}$ on
$\CC^{n+1}-\{0\}$ by $ \lambda\cdot(y_{0}, \ldots
,y_{n}):=(\lambda^{w_{0}}y_{0}, \ldots ,\lambda^{w_{n}}y_{n}).$
We denote $\PP(w)$ the quotient stack
$[\CC^{n+1}-\{0\}/\CC^{\star}]$. This stack is a
smooth proper  Deligne-Mumford stack.

For any subset $I:=\{i_{1}, \ldots ,i_{k}\}\subset \{0, \ldots ,n\}$,
we denote $w_{I}:=(w_{i_{1}}, \ldots ,w_{i_{k}})$.
We have a closed embedding  $\iota_{I}:\PP(w_{I}):=\PP(w_{i_{1}}, \ldots
,w_{i_{k}}) \to \PP(w)$. We denote $\PP(w)_{I}$ the image of this
stack morphism. In the following, we will identify $\PP(w_{I})$ with
$\PP(w)_{I}$.

Let us define the invertible sheaf $\mathcal{O}_{\PP(w)}(1)$ on
$\PP(w)$.
For any scheme $X$ and for any stack morphism $X\to \PP(w)$ given by a
principal $\CC^{\star}$-bundle $P\to X$ and a
$\CC^{\star}$-equivariant morphism $P\to \CC^{n+1}-\{0\}$, we put
$\mathcal{O}_{\PP(w)}(1)_{X}$ the sheaf of sections of the associated
line bundle of $P$.

Let us consider the following map
\begin{align*}
  \widetilde{f}_{w} : \CC^{n+1}-\{0\} &\to \CC^{n+1}-\{0\}\\
 (z_{0}, \ldots ,z_{n}) &\mapsto (z_{0}^{w_{0}}, \ldots ,z_{n}^{w_{n}})
\end{align*}
If we consider the standard action (\ie with weights $1$) on the
source of and the action with weights on the target, the map
$\widetilde{f}_{w}$ is $\CC^{\star}$-equivariant. This induces a stack
morphism $f_{w}:\PP^{n}\to \PP(w)$.  By remark $(12.5.1)$ of
\cite{LMBca}, the invertible sheaf $f^{\ast}\mathcal{O}_{\PP(w)}(1)$
is the sheaf $\mathcal{O}_{\PP^{n}}(1)$.

\subsubsection{Weighted projective spaces as orbifolds}
\label{subsubsec:weight-proj-spac}
In this part, we are using the language of orbifold used by Satake
\cite{Sgm} and W.\k Chen and Y.\k Ruan \cite{CRogw}. In this setting,
the author didn't find in the literature a complete reference for the
orbifold structure on weighted projective spaces. The purpose of this
part is to fix it in this language : namely we use the notion of good
map which is defined in \cite{CRnco}. 

First, we recall some general definitions about orbifold charts.  Let
$U$ be a connected topological space. A \emph{chart} of $U$ is a
triple $(\widetilde{U},G,\pi)$ where $\widetilde{U}$ is a connected
open set of $\CC^{n}$, $G$ is a \emph{finite} commutative \footnote{In
  the general case, one doesn't suppose that the groups are
  commutative (see \cite{CRogw}). Nevertheless, here we consider only
  examples where the groups are commutative.} group which acts
holomorphically on $\widetilde{U}$ and $\pi$ is a map from
$\widetilde{U}$ on $U$ such that $\pi$ is inducing a homeomorphism
between $\widetilde{U}/G$ and $U$. We denote $\ker(G)$ the subgroup of
$G$ that acts trivially on $\widetilde{U}$. When we will not need to
specify the group or the projection, we will denote $\widetilde{U}$
for a chart of $U$.

Let $U$ be a connected open set of $U'$. Let
$(\widetilde{U}',G',\pi')$ be a chart of $U'$. A chart $(\widetilde{U
},G,\pi)$ of $U$ is \emph{induced} by $(\widetilde{U}',G',\pi')$ if
there exists a monomorphism of groups $\kappa: G \rightarrow G'$ and
an open $\kappa$-equivariant embedding  $\alpha$ from $\widetilde{U}$
to $\widetilde{U}'$ such that $\kappa$ induces an isomorphism between
$\ker (G)$ and $\ker (G')$ and $\pi'=\alpha\circ\pi$. In \cite{Sgb},
Satake calls such pair
$(\alpha,\kappa):(\widetilde{U},G,\pi)\hookrightarrow(\widetilde{U}',G',\pi')$
an \emph{injection} of charts.

We define the action of the multiplicative group $\CC^{\star}$ on
$\CC^{n+1}-\{0\}$ by $ \lambda\cdot(y_{0}, \ldots
,y_{n}):=(\lambda^{w_{0}}y_{0}, \ldots ,\lambda^{w_{n}}y_{n}).$
The \emph{weighted projective space} is the quotient of
$\CC^{n+1}-\{0\}$ by this action. Denote by $|\PP(w)|$ this
topological space and $\pi_{w}:\CC^{n+1}-\{0\}\to |\PP(w)|$ the
quotient map. Denote by $[y_{0}:\ldots:y_{n}]$ the class of
$\pi_{w}(y_{0},\ldots ,y_{n})$ in $|\PP(w)|$. 
We have the following commutative diagram :
\begin{align*} 
  \xymatrix{(z_{0}, \ldots ,z_{n}) \ar @{|->}[d] & \CC^{n+1}-\{0\}
    \ar[d]_-{\widetilde{f}_{w}} \ar[r]^-{\pi} & \PP^{n}\ar[d]^-{f_{w}}
    & [z_{0}: \ldots :z_{n}]\ar @{|->}[d] \\ (z_{0}^{w_{0}}, \ldots
    ,z_{n}^{w_{n}})& \CC^{n+1}-\{0\} \ar[r]^-{\pi_{w}} & |\PP(w)|&
    [z_{0}^{w_{0}}: \ldots :z_{n}^{w_{n}}]}
\end{align*} 
where $\pi$ is the standard quotient map for complex projective space.
Denote by $\bs{\mu}_{k}$ the group of $k$-th roots of unity.  We can
endow $|\PP(w)|$ with two different orbifold structures.
 In the algebraic settings, we say that the
  Deligne-Mumford stacks $\PP(w)$ and
  $[\PP^{n}/\bs{\mu}_{w_{0}}\times\cdots\times\bs{\mu}_{w_{n}}]$ have
  the same coarse moduli space $|\PP(w)|$. 
\begin{enumerate}  
\item[(i)]  
   The group $\bs{\mu}_{w_{0}}\times\cdots\times\bs{\mu}_{w_{n}}$ acts 
  on $\PP^n$ in the following way : 
  \begin{align*} 
\bs{\mu}_{w_{0}}\times\cdots\times\bs{\mu}_{w_{n}}\times \PP^n &
\longrightarrow \PP^n
\\ \left((\lambda_{0},\ldots,\lambda_{n}),[z_{0}:\ldots:z_{n}]\right)&\longmapsto 
[\lambda_{0}z_{0}:\ldots:\lambda_{n}z_{n}]
\end{align*} 
 The map $f_{w}:\PP^{n}\to\PP(w)$ induces  a homeomorphism between 
$\PP^{n}/\bs{\mu}_{w_{0}}\times\cdots\times\bs{\mu}_{w_{n}}$ and
$|\PP(w)|$. So, the topological space 
$|\PP(w)|$ is endowed with an orbifold structure. 
\item[(ii)] The topological space $|\PP(w)|$ can be also endowed with
  an orbifold structure, which is not global, via the map $\pi_{w}$.
  The orbifold atlas which defines this structure is described below.
 \end{enumerate} 

 In this article, we will study only the orbifold structure which
 comes from (ii).\footnote{The author has found in the literature a
    mixing between these two orbifold structures on the same
   topological space $|\PP(w)|$. This is one of the motivation to
   make explicit the orbifold structure which comes from (ii).}  For $i \in \{0,\ldots,n\}$, denote
 $U_{i}:=\{[y_{0}:\ldots:y_{n}] \mid y_{i}\neq 0\}\subset |\PP(w)|$.
 Let $\widetilde{U}_{i}$ be the set of points of $\CC^{n+1}-\{0\}$
 such that $y_{i}=1$.  The subgroup of $\CC^\star$ which stabilizes
 $\widetilde{U}_{i}$ is $\bs{\mu}_{w_{i}}$.  The map
 $\pi_{i}:=\pi_{w}\mid_{\widetilde{U}_{i}}: \widetilde{U}_{i}
 \longrightarrow U_{i}$ induces a homeomorphism between
 $\widetilde{U}_{i}/\bs{\mu}_{w_{i}}$ and $U_{i}$.
  
 Let $U$ be a connected open set of $|\PP(w)|$. A chart
 $(\widetilde{U},G_{\widetilde{U}},\pi_{\widetilde{U}})$ of $U$ is
 called \emph{admissible} if there exists $i\in\{0, \ldots ,n\}$ such
 that $\widetilde{U}\subset \widetilde{U}_{i}$ is a connected
 component of $\pi_{i}^{-1}(U)$, $G_{\widetilde{U}}$ is the subgroup
 of $\bs{\mu}_{w_{i}}$ which stabilizes $\widetilde{U}$ and
 $\pi_{\widetilde{U}}= \pi_{i}\mid_{\widetilde{U}_{i}}$.
In particular, the charts
$(\widetilde{U_{i}},\bs{\mu}_{w_{i}},\pi_{i})$ of $U_{i}$ are
admissible charts.  Denote by $\mathcal{A}(|\PP(w)|)$ the set of all
admissible charts. The set of charts of $\mathcal{A}(|\PP(w)|)$ induces
a cover, denoted by $\mathcal{U}_{w}$, of $|\PP(w)|$.

\begin{prop}\label{prop:atlas}
  The set $\mathcal{A}(|\PP(w)|)$ is an orbifold atlas.
\end{prop} 

We will denote by $\PP(w)$ the orbifold
$(|\PP(w)|,\mathcal{A}(|\PP(w)|))$.

\begin{proof}[Proof of Proposition \ref{prop:atlas}]
  According to \cite{MPos}, we have to prove that the cover
  $\mathcal{U}_{w}$ satisfies the following conditions :
  \begin{enumerate}
  \item each open set $U$ of the cover $\mathcal{U}_{w}$ has a chart
    $(\widetilde{U},G_{U},\pi_{U})$,
  \item for any $p$ in $U\cap V$, there exists $W\subset U\cap V$ which
    contains $p$ and two injections of charts $\widetilde{W}\hookrightarrow
    \widetilde{U}$, $\widetilde{W}\hookrightarrow\widetilde{V}$.
  \end{enumerate}
  The first point is clear. Let $(\widetilde{U},G_{U},\pi_{U})$ be a
  chart of $U$ and $(\widetilde{V},G_{V},\pi_{V})$ be a chart of $V$
  in $\mathcal{A}(|\PP(w)|)$. Let $p$ be a point in $U\cap V$.  By
  definition of $\mathcal{A}(|\PP(w)|)$ there exists a unique pair
  $(i,j)\in\{0, \ldots ,n\}$ such that $\widetilde{U}\subset
  \widetilde{U}_{i}$ and $\widetilde{V}\subset \widetilde{U}_{j}$.  We
  can find a chart $(\widetilde{U}_{p},G_{p},\pi_{p})$ of a small
  neighborhood $U_{p}$ of $p$ such that $\widetilde{U}_{p}\subset
  \widetilde{U}\subset \widetilde{U}_{i}$ and the map
 \begin{align*}
  \psi_{ij}: \widetilde{U}_{p} & \longrightarrow \widetilde{V}\subset
  \widetilde{U}_{j} \\
  (y_{0}, \ldots ,1_{i}, \ldots ,y_{n})&\longmapsto
  (y_{0}/y_{j}^{w_{0}/w_{j}},\ldots,1_{j},\ldots,y_{n}/y_{j}^{w_{n}/w_{j}})
\end{align*}
where $y_{j}^{1/w_{j}}$ is a $w_{j}$-th roots of $y_{j}$, is an
injection of charts. For more details about the existence of such a
chart, see Proposition IV$.1.10$ of \cite{Mthesis}.
\end{proof}

\begin{rem}\label{rem:group,trivial}
  On an orbifold, one can define a group which acts globally and
  trivially (\cf Part $4.1$ of \cite{CRogw} or Lemma $3.1.10$ of 
\cite{Mthesis}).
For $\PP(w)$, it is easy to see that this group is $\bs{\mu}_{\pgcd(w)}$.
\end{rem}

Between two orbifolds, one can define orbifold maps (see Paragraph
$4.1$ in \cite{CRogw}). But then one has some problems when you want
to pull back bundles. So, one defines a more restrictive map which is
called \emph{good} (see Section $4.4$ of \cite{CRogw}), that allows 
to pull back bundles.

\begin{prop}\label{prop:injection,bonne}
Let $I:=\{i_{1}, \ldots ,i_{k}\}\subset \{0, \ldots ,n\}$. The
inclusion map
 \begin{align*} 
   \iota_{I} : \PP(w_{I})& \longrightarrow
   \PP(w)\\
   [z_{1}:\ldots:z_{\delta}]& \longmapsto
   [0:\ldots:0:z_{i_{1}}:0:\ldots:0:z_{i_{\delta}}:0:\ldots:0]
\end{align*} 
is a good orbifold map.
\end{prop} 

\begin{proof}
  We will prove this proposition for the set $I=\{0, \ldots
  ,\delta\}$. First we use  Section $4.1$ of \cite{CRogw} to
  construct a \emph{compatible cover} (\cf Section $4.1$ of
  \cite{CRogw}), denoted by $\mathcal{U}_{I}$, associated to the atlas
  $\mathcal{A}(|\PP(w_{I})|)$. To have a good map, we need a
  correspondence between open sets and injections of charts which
  satisfies some conditions.  We denote this correspondence by
  $\mathfrak{F}$.  For any open set $U_{I}$ of $\mathcal{U}_{I}$, we
  put $\mathfrak{F}(U_{I}):=\{[y_{0}:\ldots:y_{n}]\mid
  [y_{0}:\ldots:y_{n}]\in U_{I}\}$.  For any injection
  $(\alpha,\kappa):\widetilde{U}_{I}\hookrightarrow
  \widetilde{V}_{I}$, we put
  $\mathfrak{F}(\alpha,\kappa):=(\alpha,\id):\widetilde{\mathfrak{F}(U_{I})}\hookrightarrow
  \widetilde{\mathfrak{F}({V}_{I})}$.  It is straightforward to check
  that these data satisfy the conditions to be a good map (see
  Proposition IV$.1.15$ of \cite{Mthesis}).
\end{proof}

For any subset $I$ of $\{0, \ldots ,n\}$, we define the topological
space $|\PP(w)_{I}|:=\iota_{I}(\PP(w_{I}))$. The
orbifold atlas of $|\PP(w)|$ induces a natural orbifold atlas
(\cf Remark IV.$1.11.(5)$ of \cite{Mthesis}) which endowed
$|\PP(w)_{I}|$ by an orbifold structure denoted by $\PP(w)_{I}$.
The orbifold map
  $\iota_{I}:\PP(w_{I})\rightarrow \PP(w)$ induces an isomorphism
  between $\PP(w_{I})$ and $\PP(w)_{I}$.
In the following, we will identify the orbifolds
$\PP(w_{I})\hookrightarrow \PP(w)$ and $\PP(w)_{I}\subset \PP(w)$.
 
\begin{prop}\label{prop:f,bonne}
  The map  
  \begin{align*} 
    f_{w} \,:  \PP^{n} & \longrightarrow  \PP(w) \\
    [z_{0}: \ldots:z_{n}]  &\longmapsto  [z_{0}^{w_{0}}: \ldots:z_{n}^{w_{n}}]\\
\end{align*} 
is  a good orbifold map.
\end{prop}    

\begin{rem}\label{rem:degre} 
  The degree of $f_{w}$, regarded as a map between topological spaces, is
  $\prod w_{i}/ \gcd(w_{0}, \ldots ,w_{n})$.
\end{rem} 

\begin{proof}[Proof of Proposition \ref{prop:f,bonne}] 
  We will just prove this proposition when the weights are relatively
  prime. We refer to  Proposition IV.$1.18$ of \cite{Mthesis} for
  the general case.  Recall that the map
  \begin{align*} 
    \widetilde{f}_{w}: \CC^{n+1}-\{0\} & \longrightarrow   \CC^{n+1}-\{0\} \\
    (z_{0}, \ldots ,z_{n}) & \longmapsto (z_{0}^{w_{0}}, \ldots
    ,z_{n}^{w_{n}})
  \end{align*} is $\CC^{\star}$-equivariant and  it lifts the map 
  $|f_{w}|:|\PP^{n}|\rightarrow |\PP(w)|$.  The map $f_{w}$ is
  surjective and open.  So, the map $f_{w}$ is a regular orbifold map,
  that is $f_{w}^{-1}(|\PP(w)_{\reg}|)$ is open, connected and
  dense, where $|\PP(w)_{\reg}|$ is $\{p \in |\PP(w)|\mid G_{p}=\bs{\mu}_{\gcd(w)}\}$.
 Then, Lemma $4.4.11$ of \cite{CRogw} shows that
  the map $f_{w}$ is good.
  \end{proof}
 
\begin{prop}\label{prop:fibre,vectoriel}
  There exists an orbibundle of rank $1$, denoted by
  $\mathcal{O}_{\PP(w)}(1)$, over $\PP(w)$ such that
  $f_{w}^{\ast}\mathcal{O}_{\PP(w)}(1)$ is isomorphic to the
  bundle $\mathcal{O}_{\PP^{n}}(1)$ over $\PP^{n}$.
\end{prop} 

\begin{rem}\label{rem:orbibunble}
 For any subset  $I$ of $\{0,\ldots,n\}$, the orbibundle
 $\iota_{I}^{\ast}\mathcal{O}_{\PP(w)}(1)$ is isomorphic to 
 the orbibundle  $\mathcal{O}_{_{\PP(w_{I})}}(1)$. 
\end{rem}
                

  \begin{proof}[Proof of Proposition \ref{prop:fibre,vectoriel}]  
    We will define the orbibundle $\mathcal{O}_{\PP(w)}(1)$ on
    $\PP(w)$ by its transition maps.  Let
    $\alpha:\widetilde{U}\hookrightarrow \widetilde{V}$ be an
    injection between two charts of $\mathcal{A}(|\PP(w)|)$.  By
    definition of $\mathcal{A}(|\PP(w)|)$, there exists a unique pair
    $(i,j)\in\{0, \ldots ,n\}$ such that $\widetilde{U}\subset
    \widetilde{U}_{i}$ and $\widetilde{V}\subset \widetilde{U}_{j}$.
    For any $y\in \widetilde{U}$ and any $t\in \CC$, we put
\begin{align}\label{eq:transition,O(d)}
  \psi^{\mathcal{O}_{_{\PP(w)}}(1)}_{\alpha}(y)(t)=
  \begin{cases}
    \zeta t & \mbox{ if } i=j, \\
    t/y_{j}^{1/w_{j}} & \mbox{ if } i\neq j,
  \end{cases}
\end{align}
where $\zeta\in\bs{\mu}_{w_{i}}$.  The cocycle condition is easy to
check. So, we have defined an orbibundle of rank $1$, denoted by
$\mathcal{O}_{\PP(w)}(1)$, on $\PP(w)$.

The map $f_{w}$ is a good map (see Proposition \ref{prop:f,bonne}), so
we can define the pull back bundle
$f_{w}^{\ast}\mathcal{O}_{\PP(w)}(1)$ over $\PP^{n}$. A careful
computation shows that the bundles
$f_{w}^{\ast}\mathcal{O}_{\PP(w)}(1)\to \PP^{n}$ and
$\mathcal{O}_{\PP^{n}}(1)\to \PP^{n}$ have the same transition
maps.
\end{proof}
\subsection{Orbifold cohomology of weighted projective spaces as
  $\CC$-vector space}
\label{subsec:Orbif-cohom-weight}
We refer to  Definition $3.2.3$ of \cite{CRnco} for the definition
of orbifold cohomology. 

 For any $g\in \cup \bs{\mu}_{w_{i}}$, there exists a unique $\gamma(g)$
 in $[0,1[$ such that $g=\exp(2i\pi\gamma(g))$. When there will be no
 confusion, we will simply write $\gamma$ instead of $\gamma(g)$.

 \begin{prop}\label{prop:coho,orbifold}
   For any $g\in \cup \bs{\mu}_{w_{i}}$, we put $\age(g):=\{\gamma
   w_{0}\}+\cdots+\{\gamma w_{n}\}$ and $I(g):=\{i\mid g\in
   \bs{\mu}_{w_{i}}\}$ where $\{\cdot\}$ means the fractional part.
 The graded $\CC$-vector space structure of the
   orbifold cohomology of $\PP(w)$ is given by the following :
  \begin{align*} H^{2\star}_{\orb}(\PP(w),\CC) & =  
   \hspace{-0.3cm} \bigoplus_{g\in\cup \bs{\mu}_{w_{i}}} 
    \hspace{-0.3cm}  H^{2(\star-\age(g))} ( |\PP(w)_{(g)}|,\CC) \\ 
&\simeq \hspace{-0.3cm}\bigoplus_{g\in\cup \bs{\mu}_{w_{i}}} 
     \hspace{-0.3cm} H^{2(\star-\age(g))} ( |\PP(w_{I(g)})|,\CC) \\ 
\end{align*} 
\end{prop} 

\begin{rem}
  In \cite{Kcwps}, T.\ Kawasaki shows the following
  results
    \begin{align*} 
    H^{2i}(|\PP(w)|,\CC)=\begin{cases} \CC & \mbox{ if } 
        i\in\{0, \ldots ,n\} \,; \\ 0 & \mbox{otherwise.}  \end{cases} 
    \end{align*} 
    Then, according to the proposition above, we have an 
    explicit description of the  $\CC$-vector space  
    $H^{2\star}_{\orb}(\PP(w),\CC)$. 
\end{rem}

\begin{proof}[Proof of Proposition \ref{prop:coho,orbifold}]
  Let $p$ be in $|\PP(w)|$. Let $(\widetilde{U}_{p},G_{p},\pi_{p})$ be
  a chart of a neighborhood of $U_{p}$ of $p$. Denote by
  $\widetilde{p}$ the lift of $p$ in $\widetilde{U}_{p}$. The action
  of $G_{p}$ on the tangent space $T_{\widetilde{p}}\widetilde{U}_{p}$
  induces the following representation of the group $G_{p}$
 \begin{align*} G_{p} \longrightarrow & GL(n,\CC) \\ 
   g=\exp(2i\pi\gamma) \longmapsto & \diag (e^{2i\pi\gamma w_{0}},
   \ldots ,e^{2i\pi\gamma w_{n}})
\end{align*} 
According to  Part $3.2$ of \cite{CRnco}, the \emph{age} of
$g$ is $\{\gamma w_{0}\}+\cdots+\{\gamma w_{n}\}$.  We deduce the
equality of proposition.

For the second part, we remark that the twisted sector 
$|\PP(w)_{(g)}|$ is $|\PP(w)_{I(g)}|$ which is
identified with $|\PP(w_{I(g)})|$.
\end{proof}

For any $g$ in $\cup \bs{\mu}_{w_{i}}$, put
\begin{align*} 
\eta_{g}^{d} := 
c_{1}(\mathcal{O}_{\PP(w)_{(g)}}(1))^{d} 
\in H^{2d}(|\PP(w)_{(g)}|,\CC).\end{align*} 
Remark that 
$\eta_{g}^{d}$ vanishes for $d>\dim_{\CC}\PP(w)_{(g)}$.
We deduce the following corollary.

\begin{cor} \label{cor:degre,coho,orbifold}
  \begin{enumerate}
\item\label{item:cor:dimension,coho,orbifold} The dimension of the
  $\CC$-vector space  $H^{2\star}_{\orb}(\PP(w),\CC)$ is
  $\mu:=w_{0}+\cdots+w_{n}$. 
\item  The set  $\bs{\eta} :=\{\eta^{d}_{g}\mid g\in \cup \bs{\mu}_{w_{i}}, 
  d\in\{0, \ldots ,\dim_{\CC}\PP(w)_{(g)}\}\}$ is a basis of the
  $\CC$-vector space $H^{\star}_{\orb}(\PP(w),\CC)$. The orbifold
  degree of  
  $\eta^{d}_{g}$ is $2(d + \age(g))$. 
  \end{enumerate}
\end{cor} 

We refer to Formula $(2.4)$ of  \cite{CRnco} for the definition of the orbifold
integral (see also Formula (III$.3.3$) of \cite{Mthesis}). 

\begin{prop}\label{prop:integrale,orbifold}
 We have the following equality
  \begin{align*} 
  \int^{\orb}_{\PP(w)}\eta_{1}^{n}=\prod_{i=0}^{n} w_{i}^{-1}. 
  \end{align*} 
\end{prop} 

\begin{proof} We denote $|\PP(w)_{\reg}|:=\{p\in|\PP(w)\mid
  G_{p}=\bs{\mu}_{\gcd(w)}\}$. Let us note that $|\PP(w)_{\reg}|$ is
  open and dense in $|\PP(w)|$.
By definition of the orbifold integral and the
  Proposition \ref{prop:fibre,vectoriel}, we have
   \begin{align*} 
  \ds{\int^{\orb}_{\PP(w)}\hspace{-.4cm}\eta_{1}^{n}=\frac{1}{\pgcd(w)}
\int_{|\PP(w)_{\reg}|}\hspace{-1cm}c_{1}(\mathcal{O}_{\PP(w)}(1))^{n}=\frac{1}{\pgcd(w)\deg(f_{w})}\int_{\PP^{n}}\hspace{-.2cm}c_{1}(\mathcal{O}_{\PP^{n}}(1))^{n}.} 
  \end{align*} 
Then  the  Remark
\ref{rem:degre} implies the proposition.
\end{proof}

According to Section $3.3$ of \cite{CRnco}, we can define an orbifold
Poincar\'e duality, denoted by $\langle \cdot,\cdot \rangle$, on
orbifold cohomology.   The proposition below is a straightforward consequence of
the definition and of  Proposition \ref{prop:integrale,orbifold}.

\begin{prop}\label{prop:dualite,poincare}
  Let $\eta^{d}_{g}$ and $\eta^{d'}_{g}$ be two elements of
  the basis  $\bs{\eta}$. 
\begin{enumerate} \item  If  $g'\neq g^{-1}$ then we have 
 $\langle \eta^{d}_{g}, \eta^{d'}_{g'} \rangle = 0$. 
   
\item If $g'=g^{-1}$ then we have $I(g)=I(g')$ and 
  \begin{align*} 
  \langle \eta^{d}_{g},\eta^{d'}_{g^{-1}}\rangle= 
  \begin{cases}\ds{\prod_{i\in I(g)}w_{i}^{-1}} & \mbox{if 
      } \deg( 
    \eta^{d}_{g})+\deg(\eta^{d'}_{g^{-1}})=2n \\0& 
    \mbox{otherwise.}  \end{cases} 
  \end{align*} 
  \end{enumerate} 
\end{prop} 

\subsection{Orbifold cohomology ring of weighted projective spaces}
\label{subsec:Orbif-cohom-ring}
Before computing the orbifold cup product in the basis $\bs{\eta}$,
we will state a lemma about the orbifold tangent bundle of 
$\PP(w)$.
 
According to Proposition  \ref{prop:injection,bonne}, an orbibundle on
$\PP(w)$ can be restricted to $\PP(w_{I})$.

\begin{lem}\label{lem:decomposition,tangent}
  For any subset $I$ of $\{0, \ldots ,n\}$, we have the following decomposition : 
  \begin{align*} 
    T\PP(w)\mid_{\PP(w_{I})}\simeq \left(\bigoplus_{i\in
        I^{c}}\mathcal{O}_{\PP(w_{I})}(w_{i})\right)\bigoplus
    T\PP(w_{I})
  \end{align*} 
where $I^{c}$ is the complement of $I$ in $\{0, \ldots ,n\}$.
\end{lem} 
 
\begin{proof}
A straightforward computation shows that these two orbibundles have
  the same transition functions.  
\end{proof}

In order to compute the orbifold cup product on
$H^{\star}_{\orb}(\PP(w),\CC)$, we will compute the trilinear form
$(\cdot,\cdot,\cdot)$ introduced in \cite{CRnco} and get the cup
product through the Poincar\'e pairing by the formula
\begin{align}\label{eq:defi,cup,product}
  (\alpha_{1},\alpha_{2},\alpha_{3})&=\langle \alpha_{1}\cup\alpha_{2},\alpha_{3}\rangle
  \ \ \forall \alpha_{1},\alpha_{2},\alpha_{3} \in H^{\star}_{\orb}(\PP(w),\CC).
\end{align}

Let us fix $g_{0},g_{1},g_{\infty}$  in $\cup \bs{\mu}_{w_{i}}$ satisfying
$g_{0}g_{1}g_{\infty}=1$. 
Let us fix a presentation of
the fundamental group $\pi_{1}(\PP^{1}-\{0,1,\infty\},\star)$ as
$\langle\lambda_{0},\lambda_{1},\lambda_{\infty}\mid
\lambda_{0}\lambda_{1}\lambda_{\infty}=1\rangle$. The group
homomorphism $\pi_{1}(\PP^{1}-\{0,1,\infty\},\star)\to\CC^{\star}$
sending $\lambda_{j}$ to $g_{j}$ for any $j$ in
$\{0,1,\infty\}$, whose image is denoted by $H$, defines a covering
$\Sigma^{\circ}$ of $\PP^{1}-\{0,1,\infty\}$ having $H$ has its
automorphism group.  
This covering extends as a ramified covering $\pi:\Sigma\to\PP^{1}$,
where $\Sigma$ is a compact Riemann surface, and the action of $H$ also
extends to $\Sigma$ in such a way that $\PP^{1}=\Sigma/H$.

 

The group $H$ acts in a natural way on
$T\PP(w)\mid_{\PP(w)_{(g_{0},g_{1},g_{\infty})}}$ where
$\PP(w)_{(g_{0},g_{1},g_{\infty})}$ is the standard notation for the
triple twisted sectors.  For any $k\in\{0,1,\infty\}$, we denote by $\iota_{g_{k}}$ the  injection $\PP(w)_{(g_{k})} \hookrightarrow
\PP(w)_{(g_{0},g_{1},g_{\infty})}$. We define the following orbibundle
\begin{align*}
  E_{(g_{0},g_{1},g_{\infty})}:=\left(T\PP(w)\mid_{\PP(w)_{(g_{0},g_{1},g_{\infty})}}\otimes
  H^{0,1}(\Sigma,\CC)\right)^{H}.
\end{align*}
In the basis $\bs{\eta}$, we define the trilinear form
$(\cdot,\cdot,\cdot)$ : 
\begin{align}\label{eq:defi,3,tensor}(\eta^{d_{0}}_{g_{0}},\eta^{d_{1}}_{g_{1}},\eta^{d_{\infty}}_{g_{\infty}} 
  ) :=
  \ds{\int^{\orb}_{\PP(w)_{(g_{0},g_{1},g_{\infty})}}\hspace{-1.5cm}\iota_{g_{0}}^{\ast}\eta^{d_{0}}_{g_{0}}\wedge
    \iota_{g_{1}}^{\ast}\eta^{d_{1}}_{g_{1}}\wedge
    \iota_{g_{\infty}}^{\ast}\eta^{d_{\infty}}_{g_{\infty}}\wedge
    c_{\max}(E_{(g_{0},g_{1},g_{\infty})})}
\end{align} 

\begin{thm}\label{thm:fibre,obstruction}
  Let $g_{0},g_{1}$ and $g_{\infty}$ be in $\cup \bs{\mu}_{w_{i}}$ such
  that $g_{0}g_{1}g_{\infty}=1$.  For $i\in\{0,1,\infty\}$, we denote
  $\gamma_{i}$ the unique element in $[0,1[$ such that
  $g_{i}=\exp(2i\pi\gamma_{i})$. The orbibundle
  $E_{(g_{0},g_{1},g_{\infty})}$ is isomorphic to
   \begin{align*} 
     \bigoplus_{j\in J(g_{0},g_{1},g_{\infty})}\hspace{-.5cm}
     \mathcal{O}_{\PP(w)_{(g_{0},g_{1},g_{\infty})}}(w_{j})
  \end{align*} 
  where $J(g_{0},g_{1},g_{\infty}):=\{i\in\{0,
  \ldots ,n\}\mid
  \{\gamma_{0}w_{i}\}+\{\gamma_{1}w_{i}\}+\{\gamma_{\infty}w_{i}\}=2\}$.
\end{thm}

\begin{proof}
According to the decomposition of Lemma
  \ref{lem:decomposition,tangent}, the obstruction bundle
  $E_{(g_{0},g_{1},g_{\infty})}$ is isomorphic to
 \begin{align*}
   \left( \bigoplus_{\stackrel{i\in \{0, \ldots ,n\}-} {I(g_{0})\cap
         I(g_{1})\cap I(g_{\infty})}}
\hspace{-1cm}     \mathcal{O}_{\PP(w)_{(g_{0},g_{1},g_{\infty})}}(w_{i})\otimes
     H^{0,1}(\Sigma,\CC) \right)^{H}\hspace{-0.3cm}\bigoplus
   \Bigg(T\PP(w)_{(g_{0},g_{1},g_{\infty})} \otimes
   H^{0,1}(\Sigma,\CC)\Bigg)^{H}.
 \end{align*}
As $H$ acts trivially on
 $T\PP(w)_{(g_{0},g_{1},g_{\infty})}$, we get 
\begin{align*} 
  \ds{E_{(g_{0},g_{1},g_{\infty})}=\hspace{-1cm}\bigoplus_{\stackrel{i\in \{0, \ldots ,n\}-}
       {I(g_{0})\cap I(g_{1})\cap I(g_{\infty})}}\hspace{-.9cm}
  \left(  \mathcal{O}_{\PP(w)_{(g_{0},g_{1},g_{\infty})}}(w_{i})\otimes 
      H^{0,1}(\Sigma,\CC)\right)^{H}}. 
  \end{align*} 
  Note that the group $H$ acts on the fiber of
  $\mathcal{O}_{\PP(w)_{(g_{0},g_{1},g_{\infty})}}(w_{i})$ by
  multiplication by the character $\chi_{i}:H\to\CC^{\star}$ which sends
  $h$ to $h^{w_{i}}$. Now we apply
  the Proposition $6.3$ of \cite{BCSocdms} (see also Proposition $3.4$
  of \cite{CH} or Theorem IV.$5.13$ of \cite{Mthesis}) and we get the
  theorem.
\end{proof}

 Theorem \ref{thm:fibre,obstruction} and Formula \eqref{eq:defi,3,tensor}  of the 
trilinear form $(\cdot,\cdot,\cdot)$  imply that we can compute this trilinear form in
the basis $\bs{\eta}$. Then, the definition of the cup product via
Formula \eqref{eq:defi,cup,product} gives us the following corollary.
  
\begin{cor}\label{cor:cup,produit}
Let $\eta^{d_{0}}_{g_{0}}$ and 
  $\eta^{d_{1}}_{g_{1}}$ be two elements of the basis $\bs{\eta}$. 
  We have 
  \begin{align*} 
  \eta^{d_{0}}_{g_{0}}\cup \eta^{d_{1}}_{g_{1}}= 
  \left(\prod_{i\in K(g_{0},g_{1})} \hspace{-0.3cm}w_{i}\right) 
  \eta^{d}_{g_{0}g_{1}} 
  \end{align*} where $K(g_{0},g_{1}):= 
        J\Bigl(g_{0},g_{1},(g_{0}g_{1})^{-1}\Bigr)\bigsqcup 
        I(g_{0}g_{1})-I(g_{0})\cap 
        I(g_{1})$ 
  and 
  $d:=\frac{\deg(\eta^{d_{0}}_{g_{0}})}{2}+\frac{\deg(\eta^{d_{1}}_{g_{1}})}{2}-\age(g_{0}g_{1})=d_{0}+d_{1}+\age(g_{0})+\age(g_{1})-\age(g_{0}g_{1})$. 
\end{cor} 

 \begin{expl} 
  Let us consider the case where the weights   are
  $w=(1,2,2,3,3,3)$ (this example was considered in \cite{Jocwps}).
  The orbifold cup product is computed in the table below in  basis
  $\bs{\eta}$. We put $j:=\exp(2i\pi/3)$. 
  \begin{align*} 
    \begin{array}{|c|c|c|c|c|c|c||c|c|c|c|c|c|c|c|} 
\hline &\eta_{1}^{0}&\eta_{1}^{1}&\eta_{1}^{2}&\eta_{1}^{3}&\eta_{1}^{4}&\eta_{1}^{5}&\eta_{j}^{0}&\eta_{j}^{1}&\eta_{j}^{2}&\eta_{-1}^{0}&\eta_{-1}^{1}&\eta_{j^{2}}^{0}&\eta_{j^{2}}^{1}&\eta_{j^{2}}^{2} \\ 
\hline \eta^{0}_{1}&\eta_{1}^{0}&\eta_{1}^{1}&\eta_{1}^{2}&\eta_{1}^{3}&\eta_{1}^{4}&\eta_{1}^{5}&\eta_{j}^{0}&\eta_{j}^{1}&\eta_{j}^{2}&\eta_{-1}^{0}&\eta_{-1}^{1}&\eta_{j^{2}}^{0}&\eta_{j^{2}}^{1}&\eta_{j^{2}}^{2} \\ 
\hline\eta_{1}^{1}&&\eta_{1}^{2}&\eta_{1}^{3}&\eta_{1}^{4}&\eta_{1}^{5}&0&\eta_{j}^{1}&\eta_{j}^{2}&0&\eta_{-1}^{1}&0&\eta_{j^{2}}^{1}&\eta_{j^{2}}^{2}&0 \\ 
\hline\eta_{1}^{2}&&&\eta_{1}^{4}&\eta_{1}^{5}&0&0&\eta_{j}^{2}&0&0&0&0&\eta_{j^{2}}^{2}&0&0 \\ 
\hline\eta_{1}^{3}&&&&0&0&0&0&0&0&0&0&0&0&0 \\ 
\hline\eta_{1}^{4}&&&&&0&0&0&0&0&0&0&0&0&0 \\ 
\hline\eta_{1}^{5}&&&&&&0&0&0&0&0&0&0&0&0 \\ 
\hline\hline\eta_{j}^{0}&&&&&&&4.\eta_{j^{2}}^{2}&0&0&0&0&4\eta_{1}^{3}&4\eta_{1}^{4}&4\eta_{1}^{5} \\ 
\hline\eta_{j}^{1}&&&&&&&&0&0&0&0&4\eta_{1}^{4}&4\eta_{1}^{5}&0 \\ 
\hline\eta_{j}^{2}&&&&&&&&&0&0&0&4\eta_{1}^{5}&0&0 \\ 
\hline\eta_{-1}^{0}&&&&&&&&&&3^{3}\eta_{1}^{4}&3^{3}\eta_{1}^{5}&0&0&0 \\ 
\hline\eta_{-1}^{1}&&&&&&&&&&&0&0&0&0 \\ 
\hline\eta_{j^{2}}^{0}&&&&&&&&&&&&1.\eta_{j}^{1}&1.\eta_{j}^{2}&0 \\ 
\hline\eta_{j^{2}}^{1}&&&&&&&&&&&&&0&0\\ 
\hline\eta_{j^{2}}^{2}&&&&&&&&&&&&&&0 \\ 
 \hline   \end{array} 
  \end{align*} 
  The upper left corner is just the standard cup product on 
  $H^{\ast}(|\PP(1,2,2,3,3,3)|)$. 
\end{expl} 

\begin{expl}
  Let $\PP(w_{0},w_{1})$ be a weighted projective line. Let us denote
  by $d$ the greatest common divisor of $w_{0}$ and $w_{1}$. Choose 
  integers $m$ and $n$ such that $mw_{0}+nw_{1}=d$. Denote $g_{w_{0}}:=\exp(2i\pi
  n/w_{0})$, $g_{w_{1}}:=\exp(2i\pi m/w_{1})$ and $g_{d}:=\exp(2i\pi/d)$.
 If we put $x:=\eta_{g_{w_{0}}}^{0}$ , $y:=\eta_{g_{w_{1}}}^{0}$ and $\xi:=\eta^{0}_{g_{d}}$, we have that 
 \begin{align*}
   H^{\star}_{\orb}(\PP(w_{0},w_{1}),\CC)&=\CC[x,y,\xi]/\langle xy, w_{0}x^{w_{0}/d}-w_{1}y^{w_{1}/d}\xi^{n-m},\xi^{d}-1\rangle
 \end{align*}
This agrees completely with the computation of \cite[Section $9$]{AGVgwdms}.
\end{expl}
 
\subsection{Some initial conditions for the Frobenius manifold}
In Sections \ref{subsec:Orbif-cohom-weight} and
\ref{subsec:Orbif-cohom-ring} we have computed two initial conditions
for the Frobenius manifold namely the orbifold Poincar\'e duality
$\langle \cdot,\cdot \rangle$ and the unit $\eta_{0}^{0}$.
In this section, we will compute a third one which is
$\id-\nabla\mathfrak{E}$ where $\nabla$ is the
torsion free connection associated to the non-degenerate pairing
$\langle\cdot,\cdot\rangle$ and $\mathfrak{E}$ is the Euler field
defined in (\ref{eq:euler,field,Aside}) below.
We  start with the following lemma.
\begin{lem}\label{lem:suite,exacte}
We have the exact sequence of orbifold bundles over $\PP(w)$  
  \begin{align*} 
\xymatrix{0\ar[r]& \underline{\CC} \ar[r]^-{\Phi} & 
    \mathcal{O}_{\PP(w)}(w_{0}) \oplus \ldots \oplus 
    \mathcal{O}_{\PP(w)}(w_{n}) \ar[r]^-{\varphi}& 
    T\PP(w) \ar[r] &0} 
  \end{align*} where $\underline{\CC}$ is the orbifold trivial bundle
  of rank $1$ over $\PP(w)$. 
\end{lem} 
 
The proof is a straightforward generalization  of the proof in
Section $3$ of Chapter $3$ of \cite{GHag} (all the details are
  explained in Lemma $V.2.1$ of \cite{Mthesis}).  

Recall that $\eta_{1}^{1}=c_{1}(\mathcal{O}_{\PP(w)}(1))$.
Lemma \ref{lem:suite,exacte} and Section $4.3$ of \cite{CRogw}
(see also Proposition $\rm{III}.4.13$ in \cite{Mthesis}) imply that 
$$
c(T\PP(w)) = c(\mathcal{O}_{\PP(w)}(w_{0}) \oplus \cdots \oplus
\mathcal{O}_{\PP(w)}(w_{n}) )=
\prod_{i=0}^{n}\left(1+w_{i}\eta_{1}^{1}\right).$$ 
 We deduce that
\begin{align}\label{eq:c1,TP(w)}
  c_{1}(T\PP(w))=\mu\eta^{1}_{1}.
\end{align}
where $\mu:=w_{0}+\cdots+w_{n}$.   

For any $g\in \cup \bs{\mu}_{w_{i}}$ and any $0\leq d\leq \dim_{\CC}\PP(w)_{(g)}$, we denote
$t_{g,d}$ the coordinates of $H^{2\star}_{\orb}(\PP(w),\CC)$ with
respect to the element of the basis $\eta_{g}^{d}$.
As in the book of Y.\k Manin \cite[p.37]{Mfm}, 
we define the Euler field by 
\begin{align}
  \label{eq:euler,field,Aside}
  \mathfrak{E}:=\hspace{-1cm}\sum_{\stackrel{g\in \cup \bs{\mu}_{w_{i}}}{0\leq d\leq \dim_{\CC}\PP(w)_{(g)}}}\hspace{-.8cm}(1-\deg(\eta^{d}_{g})/2)t_{g,d}\partial_{t_{g,d}}+\mu\partial_{t_{1,1}}.
\end{align}

\begin{notation}\label{not:order}
  For the following, it will be useful to have an order on the basis
  $\bs{\eta}$ of $H^{2\star}_{\orb}(\PP(w),\CC)$. Choose any
  determination of the argument in $\CC$. Hence, for any $g\in\cup
  \bs{\mu}_{w_{i}}$, there exists a unique $\gamma(g)\in[0,1[$ such
  that $g=\exp(2i\pi\gamma(g))$. We say that $\eta_{g}^{d}\leq
  \eta_{g'}^{d'}$ if $\gamma(g)<\gamma(g')$ or $\gamma(g)=\gamma(g')$
  and $d\leq d'$.  We denote by $g_{\max}$ the greatest element in
  $\cup \bs{\mu}_{w_{i}}$.
\end{notation}

\begin{prop}\label{prop:matrice,Ainfty} Denote by $d_{\max}$ the
  complex dimension of the twisted sector $\PP(w)_{(g_{\max})}$.
The matrix $A_{\infty}:=\id-\nabla \mathfrak{E}$ in the basis $\bs{\eta}$ 
is $   \diag (\deg(\eta_{0}^{0})/2, \ldots ,\deg(\eta^{d_{max}}_{g_{max}})/2)$ where $\nabla$ is the
torsion free connection associated to the non-degenerate pairing $\langle\cdot,\cdot\rangle$.
This matrix satisfies  $A_{\infty}+A_{\infty}^{\ast}=n\id$ where
$A_{\infty}^{\ast}$ is the adjoint of $A_{\infty}$ with respect to the
non-degenerate bilinear form $\langle\cdot,\cdot\rangle$. 
\end{prop}

\begin{proof}
   By definition, we have 
    \begin{align*}
      A_{\infty}=\frac{1}{2} \diag \left(\deg\left(\eta_{0}^{0}\right), \ldots ,\deg\left(\eta^{n}_{0}\right),\deg\left(\eta^{0}_{g_{\max}^{-1}}\right),\ldots,\deg\left(\eta^{d_{\max}}_{g_{\max}}\right)\right).
    \end{align*}
The adjoint matrix of
    $A_{\infty}$ with respect to the non degenerate bilinear form
    $\langle\cdot,\cdot\rangle$ is 
    \begin{align*}
      A_{\infty}^{\ast}= \frac{1}{2}\diag \left(\deg\left(\eta_{0}^{n}\right), \ldots ,\deg\left(\eta^{0}_{0}\right),\deg\left(\eta^{d_{\max}}_{g_{\max}}\right),\ldots,\deg\left(\eta^{0}_{g_{max}^{-1}}\right)\right)
    \end{align*} In order to end the proof, it is enough to check that 
    $$\deg(\eta^{d}_{g})+\deg\left(\eta^{\dim_{\CC}\PP(w)_{(g)}-d}_{g^{-1}}\right)=2n.$$
 \end{proof}

\section{Orbifold quantum cohomology of weighted projective spaces}
\label{sec:Orbif-quant-cohom}
The orbifold cohomology algebra, with its Poincar\'e pairing, of
weighted projective space is now completely determined.  We have
computed three out of  four initial conditions of the Frobenius manifold.
In this section, we will study the last initial
condition $\mathfrak{E}\star\mid_{\bs{t}=\bs{0}}$.

\subsection{The orbifold  Gromov-Witten invariants}
\label{subsec:orbif-Grom-Witt}

 First, we recall the definition of orbifold stable maps to $\PP(w)$ (see
 Paragraph $2.3$ of \cite{CRogw} for details).
 
 A good orbifold map $f$ between the orbifolds $X$ and $Y$ is an
 orbifold map and a compatible structure. A compatible structure is a
 correspondence between open sets and injections of charts of $X$ and
 open sets and injections of charts of $Y$ which satisfies some
 conditions (see Section $4.4$ of \cite{CRogw} for more details). In
 particular, the compatible structure induces a homomorphism between
 the local groups that is for any $x\in X$ we have a morphism of group
 $G_{x}\to G_{f(x)}$. We will not be more precise because we will not
 use explicitly this notion.

\begin{defi}\label{defi:stable,map}
   \emph{An orbifold stable  map to $\PP(w)$} consists of  the following data  
  \begin{itemize} 
  \item a nodal orbicurve  $(C,\bs{z},\bs{m},\bs{n})$ where
    $\bs{z}:=(z_{1}, \ldots ,z_{k})$ are $k$ distinct marked points
    such that $G_{z_{i}}=\bs{\mu}_{m_{i}}$ for $i\in\{1, \ldots ,k\}$
    and the $j$-th nodal point has the action of $\bs{\mu}_{n_{j}}$
    ;  
  \item a continuous map $f:C\rightarrow \PP(w)$   
  \item and an isomorphism class of compatible structure, denoted by
    $\xi$.
  \end{itemize} 
These data  $(f,(C,\bs{z},\bs{m},\bs{n}),\xi)$ satisfy  
\begin{enumerate} 
\item for each $i$ in $I$, the orbifold map
  $f_{i}:=f\circ\varphi_{i}:C_{i}\rightarrow \PP(w)$ is holomorphic ;  
\item  for each marked or nodal point $z_{i}$, the morphism of group
  induced by $\xi$ from $G_{z_{i}}$ to 
  $G_{f(z_{i})}$ is injective 
\item and if the map $f_{i}:C_{i}\rightarrow \PP(w)$ is constant then the curve $C_{i}$ 
  has more than three singular points (i.e. nodal or marked). 
\end{enumerate} 
\end{defi} 
 
We endow the set of orbifold stable maps with the standard
equivalence relation.
Denote by $[f,(C,\bs{z},\bs{m},\bs{n}),\xi]$ the equivalence class of
the orbifold stable map $(f,(C, \bs{z},\bs{m},\bs{n}),\xi)$.
 
Let  $(f,(C,\bs{z},\bs{m},\bs{n}),\xi)$ be a stable map. We can associate 
to this stable map a homology class in 
$H_{2}(|\PP(w)|,\ZZ)$ defined by 
$f_{\ast}([C]):=\sum_{i}(f\circ\varphi_{i})_{\ast}[C_{i}]$ where
$[C_{i}]$ is the fundamental class of the curve $C_{i}$. 
This homology class does not depend on the equivalence class of the
stable map. 
For each  marked point  $z_{i}$, the  class of  compatible structure  $\xi$ induces 
a monomorphism of groups $\kappa_{i}:G_{z_{i}}\hookrightarrow 
G_{f(z_{i})}$. This  monomorphism  depends only on the equivalence
class of the stable map. 
 
Let us define the inertia orbifold  $\mathcal{I}{\PP(w)}$ by $\bigsqcup_{g\in \cup \bs{\mu}_{w_{i}}}
\PP(w)_{(g)}\times \{g\}.$ We have an evaluation map,
denoted by $\ev$, which maps a class
$[f,(C,\bs{z},\bs{m},\bs{n}),\xi]$ of stable map to
$$\left((f(z_{1}),\kappa_{1}(e^{{2i\pi}/{m_{1}}})), \ldots 
  ,(f(z_{k}),\kappa_{k}(e^{{2i\pi}/{m_{k}}}))\right) \in {\mathcal{I} \PP(w)}\times\cdots\times\mathcal{I}\PP(w).$$ 
 
An  orbifold stable map $(f,(C,\bs{z},\bs{m},\bs{n}),\xi)$ 
is said to have  \emph{type} $(g_{1}, \ldots ,g_{k})\in
\cup \bs{\mu}_{w_{i}}$ if for any $\ell\in\{1, \ldots ,k\}$,
$(f(z_{\ell}),\kappa_{1}(e^{{2i\pi}/{m_{\ell}}}))$ belongs to  
$\PP(w)_{(g_{\ell})}\times\{g_{\ell}\}$. When there is no ambiguity in
the notation, we will write  $\underline{g}$ for the $k$-uple $(g_{1}, \ldots ,g_{k})$. 
 
\begin{defi}\label{defi:moduli,space}
   Let $A$ be in $H_{2}(|\PP(w)|,\ZZ)$. We define   
$\overline{\mathcal{M}}_{k}(A,\underline{g})$ the moduli space of
equivalence classes of orbifold stable maps with $k$ marked points,
  of homology class  $A$ and of type $\underline{g}$, i.e.  
\begin{align*} 
\overline{ \mathcal{M}}_{k}(A,\underline{g})= 
\left\{  
 \begin{array}{c} 
 \mbox{ } [(f,(C,\bs{z},\bs{m},\bs{n}),\xi)]\mid 
   \# \bs{z} =k, f_{\ast}[C]=A, 
  \\ \ev(f,(C,\bs{z},\bs{m},\bs{n}),\xi)\in \prod_{\ell=1}^{k}(\PP(w)_{(g_{\ell})}\times\{g_{\ell}\})
\end{array}   
\right\}.  
\end{align*}  
 \end{defi} 

According to the results of \cite{CRogw} (\cf Proposition 
$2.3.8$), the moduli space $\overline{\mathcal{M}}_{k}(A,\underline{g})$ 
is  compact and metrizable. 
 Chen and Ruan define also a Kuranishi structure on this moduli space
 whose dimension is given by the following theorem.

 \begin{thm}[\cf Theorem $A$ of \cite{CRogw}] \label{thm:dimension,virtuelle}  
   The dimension of the Kuranishi structure described by W.\k Chen and
   Y.\k Ruan
   of $\overline{\mathcal{M}}_{k}(A,\underline{g})$ is
  \begin{align*} 
  2\left(\int_{A}c_{1}(T\PP(w))+\dim_{\CC}\PP(w)-3+k-\sum_{\ell=1}^{k}\age(g_{\ell})\right). 
  \end{align*} 
\end{thm} 
 
This Kuranishi structure defines (\cf Theorem
$6.12$ and  Section $17$ of  \cite{FOgwi}) a homology
class, called fundamental class of the Kuranishi structure,
\begin{align}\label{eq:deg,class,virtual}
   \rm{ev}_{\ast}[\overline{\mathcal{M}}_{k}(A,\underline{g})]&\in 
   H_{2\left(\int_{A}c_{1}(T\PP(w))+n-3+k-\sum_{\ell=1}^{k}\age(g_{\ell})\right)} 
 (\PP(w)_{(g_{1})}\times\cdots\times \PP(w)_{(g_{k})},\CC). 
 \end{align}

\begin{rem}\label{rem:tendentious}
   This notation is a bit tendentious because the class
  $ \rm{ev}_{\ast}[\overline{\mathcal{M}}_{k}(A,\underline{g})]$ is not
   a push forward of a homology class of
   $\overline{\mathcal{M}}_{k}(A,\underline{g})$. However, in
   algebraic geometry, we can construct a virtual fundamental class in
   the Chow group of $\overline{\mathcal{M}}_{k}(A,\underline{g})$ with
   the same degree (see Section $5.3$ in \cite[Section $4.5$]{AGVgwdms}).
 \end{rem}

For each $\ell\in\{1, \ldots ,k\}$, let $\alpha_{\ell}$ be a  class in  
$H^{2(\star-\age(g_{\ell}))}(\PP(w)_{(g_{\ell})},\CC) \subset 
H^{2\star}_{\orb}(\PP(w),\CC)$. Formula  $(1.3)$ of \cite{CRogw} 
defines the orbifold Gromov-Witten invariants by 
 \begin{align}\label{eq:defi,IGW} 
   \Psi^{A}_{k,\underline{g}} :
   H^{\star}(\PP(w)_{(g_{1})},\CC)\otimes \cdots \otimes
   H^{\star}(\PP(w)_{(g_{k})},\CC) &
   \longrightarrow  \CC \\
   \alpha_{1}\otimes\cdots\otimes\alpha_{k} & \longmapsto
   \int_{\rm{ev}_{\ast}[\overline{\mathcal{M}}_{k}(A,\underline{g})]}\alpha_{1}\wedge\cdots\wedge\alpha_{k}\nonumber.
 \end{align}


 
Let $g_{1}, \ldots ,g_{k}$ be in $\cup\bs{\mu}_{w_{i}}$. 
 Recall that $\mu:=w_{0}+\cdots+w_{n}$. Theorem \ref{thm:dimension,virtuelle} and  Formula
 (\ref{eq:c1,TP(w)}) imply that 
\begin{align*} 
  \deg \ev_{\ast}\left[\overline{\mathcal{M}}_{k}(A,\underline{g})\right]=2\left(
    \mu\int_{A}\eta_{1} +n -3 +k - \sum_{\ell=1}^{n}\age(g_{\ell}) \right).
\end{align*}

We recall that for any $g\in \cup \bs{\mu}_{w_{i}}$ and any $0\leq
d\leq \dim_{\CC}\PP(w)_{(g)}$, we denote $t_{g,d}$ the coordinate of
$H^{2\star}_{\orb}(\PP(w),\CC)$ with respect to the element of the
basis $\eta_{g,d}$. Let us put
\begin{align*}
  T:= \hspace{-1cm}\sum_{\stackrel{g\in\cup\bs{\mu}_{w_{i}}}{0\leq
      d\leq \dim_{\CC}\PP(w)_{(g)}}}\hspace{-.8cm}t_{g,d}\ \eta_{g}^{d}.
\end{align*}
The full Gromov-Witten potential of genus $0$ of the
weighted projective space $\PP(w)$, denoted by $F^{GW}$, is defined by
\begin{align*} 
F^{GW}:=\sum_{k\geq 0} \sum_{\stackrel{A\in 
    H_{2}(\PP(w),\ZZ),}{\underline{g}\in \cup (\bs{\mu}_{w_{i}})^{k}}} \hspace{-.3cm}
\frac{\Psi^{A}_{k,\underline{g}}(T, \ldots ,T)}{k!}. 
\end{align*}  


We define the orbifold quantum  product by the equation  
\begin{align}\label{eq:quantum,product} 
 \frac{\partial^3 F^{GW}(\bs{t})}{\partial 
 t_{g,d}\partial t_{g',d'} \partial t_{g'',d''}}&=\langle \partial t_{g,d}\star 
 \partial t_{g',d'}, \partial t_{g'',d''}\rangle.  
\end{align}

 \subsection{Computation of some orbifold Gromov-Witten invariants}
\label{subsec:calcul-de-certains} 


To compute the last initial condition of the Frobenius manifold, we
should compute the matrix
$A_{0}^{\circ}:=\mathfrak{E}\star\mid_{\bs{t}=0}$ (\cf Equation
(\ref{eq:quantum,product}) for the definition of the quantum orbifold
product). As the Euler field restricted to $\bs{t}=0$ is $\mu\partial
t_{1,1}$, we have to compute the Gromov-Witten invariants
$\Psi_{3}^{A}(\eta^{1}_{1},\eta^{d}_{g},\eta^{d'}_{g'})$ for any
$g,g'\in\cup \bs{\mu}_{w_{i}}$, for any $(d,d')\in \{0, \ldots
,\dim_{\CC}\PP(w)_{(g)}\}\times \{0, \ldots
,\dim_{\CC}\PP(w)_{(g')}\}$ and for any class $A \in
H_{2}(\PP(w),\ZZ)$.  By definition of Gromov-Witten invariant, if the class
$A\in H_{2}(|\PP(w)|,\ZZ)$ does not satisfy
    \begin{align}\label{eq:2}
      \mu\int_{A}\eta_{1}^{1}=
1+\deg(\eta_{g}^{d})/2+\deg(\eta_{g'}^{d'})/2-n
    \end{align}
    the Gromov-Witten invariant
    $\Psi_{3}^{A}(\eta^{1}_{1},\eta^{d}_{g},\eta^{d'}_{g'})$ is zero.
\begin{notation}\label{not:A}
In the following we denote by $A(g,d,g',d')$ the unique class in $H_{2}(|\PP(w)|,\ZZ)$
which satisfies (\ref{eq:2}). When there will be no confusion, we
denote this class by $\widetilde{A}$. 
\end{notation}
Motivated by  Corollary \ref{cor:tenseur} for the B side, we
decompose this set of Gromov-Witten invariants in the following three subsets :

\begin{align}
&\left\{ \begin{array}{l} \Psi_{3}^{\widetilde{A}}(\eta^{1}_{1},\eta^{d}_{g},\eta^{d'}_{g'}) \mbox{
    such that } \\
  1+\deg(\eta_{g}^{d})/2+\deg(\eta_{g'}^{d'})/2-n+\mu(\gamma(g^{-1})+\gamma(g'^{-1})) \neq 0 \mod
  \mu 
\end{array}
\right\} \label{eq:cond,1} \\
&\left\{ \begin{array}{l} \Psi_{3}^{\widetilde{A}}(\eta^{1}_{1},\eta^{d}_{g},\eta^{d'}_{g'}) \mbox{
    such that}  \\  1+\deg(\eta_{g}^{d})/2+\deg(\eta_{g'}^{d'})/2-n+\mu(\gamma(g^{-1})+\gamma(g'^{-1})) = 0 \mod
  \mu \\ 
 \mbox{ and } 2+\deg(\eta_{g}^{d})+\deg\eta_{g'}^{d'}=
  2n,\\   
\end{array}
\right\}\label{eq:cond,2}\\
&\left\{ \begin{array}{l} \Psi_{3}^{\widetilde{A}}(\eta^{1}_{1},\eta^{d}_{g},\eta^{d'}_{g'}) \mbox{
 such that } \\
1+\deg(\eta_{g}^{d})/2+\deg(\eta_{g'}^{d'})/2-n+\mu(\gamma(g^{-1})+\gamma(g'^{-1})) = 0 \mod
  \mu\\
\mbox{ and } 2+\deg(\eta_{g}^{d})+\deg\eta_{g'}^{d'}\neq
  2n. 
\end{array}
\right\}\label{eq:cond,3}
\end{align}

\begin{rem}\label{rem:comb}
  \begin{enumerate}
  \item\label{item:1} The  number
    $1+\deg(\eta_{g}^{d})/2+\deg(\eta_{g'}^{d'})/2-n+\mu(\gamma(g^{-1})+\gamma(g'^{-1}))$
    is equal to the integer
    \begin{align*}
1+d+d'+n-\dim\PP(w)_{(g)}-\dim\PP(w)_{(g')}+\sum_{i=0}^{n}[\gamma(g^{-1})w_{i}]+[\gamma(g'^{-1})w_{i}]      
    \end{align*} where $[\cdot]$ is the integer part.
  \item\label{item:2} Conditions (\ref{eq:cond,2}) are equivalent to the
    conditions $ gg'=1$ and $2+\deg(\eta_{g}^{d})+\deg(\eta_{g'}^{d'})=2n$.
  \end{enumerate}
\end{rem}

First we study the set (\ref{eq:cond,1}) of Gromov-Witten
invariants. The following proposition is a straightforward consequence
of Proposition \ref{prop:mult,eta,1,1}.
 
 \begin{prop}\label{prop:invariant,nul}
   Let $g,g'$ be in $\cup \bs{\mu}_{w_{i}}$ and let $(d,d')$ be in $\{0, \ldots
   ,\dim\PP(w)_{(g)}\}\times \{0, \ldots ,\dim\PP(w)_{(g')}\}$ such
   that $
   1+\deg(\eta_{g}^{d})/2+\deg(\eta_{g'}^{d'})/2-n+\mu(\gamma(g^{-1})+\gamma(g'^{-1}))
   \neq 0 \mod \mu$.  We have
   $\Psi^{\widetilde{A}}_{3}(\eta^{1}_{1},\eta^{d}_{g},\eta^{d'}_{g'})=0$
   where $\widetilde{A}$ is defined by $ \mu\int_{\widetilde{A}}\eta_{1}^{1}=
   1+\deg(\eta_{g}^{d})/2+\deg(\eta_{g'}^{d'})/2-n$.
 \end{prop}
 

 \begin{rem}
 The proposition above is proved in \cite[Theorem V.3.3 p.98]{Mthesis}
  when $\mu$ and  $\lcm(w_{0}, \ldots ,w_{n})$ are coprime. 
 \end{rem}

The following proposition computes the Gromov-Witten invariant for the
subset defined by Conditions (\ref{eq:cond,2}). According to Remark
\ref{rem:comb}.(\ref{item:2}), these conditions are
equivalent to the hypothesis of the theorem below.

\begin{prop}\label{prop:invariant,degre,0}
  Let $g,g'$ be in $\cup \bs{\mu}_{w_{i}}$ and let $(d,d')$ be in $\{0, \ldots
  ,\dim\PP(w)_{(g)}\}\times\{0, \ldots
  ,\dim\PP(w)_{(g')}\}$ such that $g'g=\id$ and
  $2+\deg\eta^{d}_{g}+\deg\eta^{d'}_{g'}=2n $. Let $A$ be a
  class in $H_{2}(|\PP(w)|,\ZZ)$.  We have
   \begin{align*} 
   \Psi^{A}_{3}(\eta^{1}_{1},\eta_{g}^{d},\eta_{g'}^{d'})=  
\begin{cases} 
\prod_{i\in I(g)}w_{i}^{-1} & 
\mbox{ if }  A=0 
\\ 
0 & \mbox{otherwise.} 
  \end{cases} 
  \end{align*} 
\end{prop} 

\begin{proof} If this invariant is not zero the class $A$ should
  satisfy the equality
  $2\mu\int_{A(g,d,g',d')}\eta^{1}_{1}=2+\deg\eta_{g}^{d}+\deg\eta_{g'}^{d'}-2n$.
By hypothesis this implies that $A=0$.
Hence   Theorem \ref{thm:fibre,obstruction} and Formula
(\ref{eq:defi,3,tensor}) finish the proof.
 \end{proof} 
 
 In order to simplify Conditions (\ref{eq:cond,3}), we will recall
 some combinatorics.  Let us denote the elements of
 $\cup\bs{\mu}_{w_{i}}$ by $1=g_{0}<g_{1}<\cdots<g_{\delta}$ where the
 order is defined by choosing the principal determination of the
 argument (\cf Notation \ref{not:order}).  Let us fix
 $g_{k}\in\cup\bs{\mu}_{w_{i}}$. There exists a unique
 triple $(d,g',d')$ in $\{0, \ldots ,\dim\PP(w)_{(g_{k})}\}\times
 \cup\bs{\mu}_{w_{i}}\times\{0, \ldots ,\dim\PP(w)_{(g')}\}$ that
 satisfies
\begin{align*}
  \begin{cases}
1+\deg(\eta_{g_{k}}^{d})/2+\deg(\eta_{g'}^{d'})/2-n+\mu(\gamma(g_{k}^{-1})+\gamma(g'^{-1}))
 = 0 \mod  \mu\\
\mbox{ and } 2+\deg\eta_{g_{k}}^{d}+\deg\eta_{g'}^{d'}\neq
  2n.  
\end{cases}
\end{align*}
Such a triple is given by
$(\dim\PP(w)_{(g_{k})},g_{k-1}^{-1},\dim\PP(w)_{(g_{k-1})})$ where $g_{-1}:=g_{\delta}$.
 
\begin{prop}\label{prop:invariants,durs}  Let $g_{k}\in\cup
  \bs{\mu}_{w_{i}}$.
  Let $\widetilde{A}$ be the class in $H_{2}(|\PP(w)|,\ZZ)$ defined by
  $\mu\int_{\widetilde{A}}\eta_{1}^{1}=
  1+\deg(\eta_{g}^{d})/2+\deg(\eta_{g'}^{d'})/2-n$. We have that
\begin{align*}
\ev_{\ast}\left[\overline{\mathcal{M}}_{2}(\widetilde{A},g_{k},g_{k-1}^{-1})\right]=(\gamma(g_{k})-\gamma(g_{k-1}))^{-1}[\PP(w)_{(g_{k})}\times\PP(w)_{(g_{k-1})}].
\end{align*}
\end{prop}

\begin{rem} \label{rem:reformulation, conj}
In \cite[Section $9$]{AGVgwdms}, D.\k Abramovich, T.\k Graber and A.\k
Vistoli have computed the small quantum cohomology of $\PP(w_{0},w_{1})$.
This result implies Proposition \ref{prop:mult,eta,1,1}, hence Propositions
\ref{prop:invariant,nul} and \ref{prop:invariants,durs} for weighted
projective lines.
\end{rem}

\begin{proof}[Proof of Proposition \ref{prop:invariants,durs}]
The divisor axiom implies that this proposition is equivalent to 
   \begin{align*} 
    \Psi^{\widetilde{A}}_{3}\left(\eta^{1}_{1},\eta^{\dim\PP(w)_{(g_{k})}}_{g_{k}},
\eta^{\dim\PP(w)_{(g_{k-1})}}_{g_{k-1}^{-1}}\right)=\hspace{-.5cm}
    \ds{\prod_{i\in I(g_{k})\bigsqcup
        I(g_{k-1})}\hspace{-.8cm}w_{i}}^{-1}.
\end{align*} 
Formula (\ref{eq:quantum,product}) and Proposition
\ref{prop:mult,eta,1,1} imply the equality above.  
\end{proof}


 Let us put
\begin{align*} 
(\!(\eta^{1}_{1},\eta^{d}_{g},\eta^{d'}_{g'})\!):=\frac{\partial^{3}F^{GW}}{\partial t_{1,1}{\partial t_{g,d}}{\partial t_{g',d'}}}\mid_{\bs{t=0}} 
\end{align*}  
Propositions \ref{prop:invariant,nul}, \ref{prop:invariant,degre,0}
and \ref{prop:invariants,durs} imply the following corollary. 
 
 \begin{cor}\label{cor:3,tenseur}Let
   $g,g'$ be in $\cup \bs{\mu}_{w_{i}}$ and let $(d,d')$ be in $\{0, \ldots
   ,\dim_{\CC}\PP(w)_{(g)}\}\times \{0, \ldots ,\dim_{\CC}\PP(w)_{(g')}\}$.
\begin{enumerate} 
\item If $
  1+\deg(\eta_{g}^{d})/2+\deg(\eta_{g'}^{d'})/2-n+\mu(\gamma(g^{-1})+\gamma(g'^{-1}))
  \neq 0 \mod \mu$ then
  $(\!(\eta^{1}_{1},\eta^{d}_{g},\eta^{d'}_{g'})\!)=0$.
\item If $
  1+\deg(\eta_{g}^{d})/2+\deg(\eta_{g'}^{d'})/2-n+\mu(\gamma(g^{-1})+\gamma(g'^{-1}))
  = 0 \mod \mu$ then we have
\begin{align*}(\!(\eta^{1}_{1},\eta^{d}_{g},\eta^{d'}_{g'})\!)= 
 \begin{cases} 
   \left({{\prod_{i\in I(g)\coprod I(g')}w_{i}}}\right)^{-1} & \mbox{ if } 
   2+\deg\eta_{g}^{d}+\deg\eta_{g'}^{d'}\neq 2n\\ 
   \left({{\prod_{i\in I(g)}w_{i}}}\right)^{-1} & \mbox{ if } 
   2+\deg\eta_{g}^{d}+\deg\eta_{g'}^{d'}= 2n
 \end{cases} 
\end{align*} 
\end{enumerate} 
\end{cor} 
 
To determine the matrix $A_{0}^{\circ}=\mathfrak{E}\star
\mid_{\bs{t}=0}$, we use Formula
(\ref{eq:quantum,product}). This Formula shows that the
data $(\!(\eta^{1}_{1},\eta^{d}_{g},\eta^{d'}_{g'})\!)$ and the
orbifold Poincar\'e duality $\langle \cdot,\cdot\rangle$ enable us to
compute the matrix $A_{0}^{\circ}$.

\section{Frobenius structure associated  to the Laurent polynomial $f$}
\label{sec:Frob-struct-assoc}
 
In this part, we will use the following notations.  
\begin{notation}\label{not:sigma}
  Let $n$ and
$w_{0},\ldots,w_{n}$ be some integers greater or equal to one.  We put
$\mu:=w_{0}+\cdots+w_{n}$.   
Consider the set $\bigsqcup_{i=0}^{n}\{\ell/w_{i} \mid
\ell\in\{0,\ldots,w_{i}-1\}\}$ where $\bigsqcup$ means the disjoint
union.   Choose a bijection $ s:\{0, \ldots ,\mu-1\}
\rightarrow \bigsqcup_{i=0}^{n}\{\ell/w_{i} \mid
\ell\in\{0,\ldots,w_{i}-1\}\} $ such that $\mathcal{V}\circ s$
is nondecreasing.  
  Let us consider, as in   \cite{DSgm2}, the 
  rational numbers $\sigma(i):=i-\mu s(i)$  for $i\in\{0,\ldots,\mu-1\}$. 
  \end{notation}
 
Let $U :=\{ (u_{0}, \ldots ,u_{n})\in \CC^{n+1}\mid u_{0}^{w_{0}}\cdots 
u_{n}^{w_{n}}=1\}$. Let  $f:U\rightarrow \CC$ be the function defined by 
$f(u_{0},\ldots,u_{n})=u_{0}+\cdots+u_{n}$.  The polyn{o}mial $f$ is not 
exactly the one considered in \cite{DSgm2} but we can apply the same methods. 
 
An easy computation shows that the critical value of $f$ are $
\mu\zeta\left({\prod_{i=0}^{n} w_{i}^{w_{i}}}\right)^{-1/\mu}$ where
$\zeta$ is a $\mu$-th roots of unity.  In \cite{DSgm1}, there
exists a Frobenius structure on the base space of any universal
unfolding of $f$.  In this example, we will see (\cf Theorem
\ref{thm:cond,init,sing}) that the Frobenius structure is semi simple
\ie we can  reconstruct the Frobenius structure  from some initial data.

 Let $A_{0}^{\circ}$ be the matrix of size $\mu\times\mu$
defined by (recall that $\overline{a+b}=a+b\mod \mu$)
\begin{align*} 
(A_{0})_{\overline{j+1},j}= 
\begin{cases}  
\mu & \mbox{ if }  s(\overline{j+1})=s(j)\,; \\ 
{\mu}\prod_{i\in I(s(j))} w_{i}^{-1} & \mbox{otherwise}. 
 \end{cases} 
\end{align*} 
The eigenvalues of $A_{0}^{\circ}$ are exactly the critical values of
$f$.  Hence, $A_{0}^{\circ}$ is a semi-simple regular matrix. In the
canonical basis $(e_{0}, \ldots ,e_{\mu-1})$ of $\CC^{\mu}$, we define
the bilinear non degenerated form $g$ by
\begin{align}\label{eq:metric} 
  g(e_{j},e_{k}) & = \begin{cases} \prod_{i\in I(s(k))}
      w_{i}^{-1} & \mbox{ if } \overline{j+k}=n\,; \\ 0 &
    \mbox{otherwise.}  \end{cases}
\end{align} 
Let $A_{\infty}$ be the matrix of size $\mu\times\mu$ defined by $
A_{\infty}=\diag (\sigma(0), \ldots ,\sigma(\mu-1))$ (\cf Section
\ref{subsec:Comb-numb-sigma} for the definition of $\sigma(\cdot)$).
This matrix satisfies $A_{\infty}+ A_{\infty}^{\ast} =n\cdot{\id}$
where $A_{\infty}^{\ast}$ is the adjoint of $A_{\infty}$ with respect
to $g$.  

\begin{thm}[Theorem $2$ of \cite{DSgm2}]\label{thm:cond,init,sing}
  The canonical Frobenius structure on any germ of a universal
  unfolding of the Laurent polynomial $f(u_{0}, \ldots
  ,u_{n})=u_{0}+\cdots+u_{n}$ on $U$ is isomorphic to the germ of
  universal semi-simple Frobenius structure with initial data
  $(A_{0}^{\circ},A_{\infty},e_{0},g)$ at the point
  \begin{align*}\left({\mu}{\prod_{i=0}^{n}w_{i}^{-w_{i}}}, 
      {\mu\zeta}{\prod_{i=0}^{n}w_{i}^{-w_{i}}},
      \ldots
      ,{\mu\zeta^{\mu-1}}{\prod_{i=0}^{n}w_{i}^{-w_{i}}}
    \right)\in\CC^{\mu}.
 \end{align*} 
\end{thm}

\subsection[The Gauss-Manin system and the Brieskorn lattice of
$f$]{The Gauss-Manin system and the Brieskorn lattice of the Laurent
  polynomial $f$}
\label{subsec:Gauss-Manin-system}

One can compute the initial data $(A_{0}^{\circ},A_{\infty},e_{0},g)$ of the Frobenius
structure from the Jacobian algebra of $f$.
Namely, the product on the Frobenius manifold at the
origin comes from the product on the Jacobian algebra of $f$ via the
isomorphism given by the primitive form
\begin{align*} 
 \omega_{0}:= \frac{\frac{d u_{0}}{u_{0}}\wedge \cdots \wedge \frac{d u_{n}}{u_{n}}}{d\left(\prod_{i}u_{i}^{w_{i}}\right)}\left|_{\prod_{i}u_{i}^{w_{i}}=1}  \right.
\end{align*} 
Moreover, the multiplication by the Euler field is induced, via this
isomorphism, by the multiplication by $f$. Finally, the non degenerated
form $g$ is given by a residue formula.
In this example, the form $g$ can also be computed  from a
duality on the Brieskorn lattice of $f$. In this example, this way is
easier. For this reason, we will use the Gauss-Manin system and the
Brieskorn lattice of $f$ to get the initial data.

 In the following, we will not give details, we refer to
 \cite{DSgm1}, \cite{DSgm2}. 
For a better exposition, we will suppose that the weights are relatively
prime (see \cite{Mthesis} for the general case).
The Gauss-Manin system of $f$ is defined by $ G:=
\Omega^{n}(U)[\theta,\theta^{-1}]/(\theta d - d
f\wedge)\Omega^{n-1}(U)[\theta,\theta^{-1}].$
The Brieskorn lattice of $f$, defined by $G_{0} :=
Im(\Omega^{n}(U)[\theta]\rightarrow G)$, is a free
$\CC[\theta]$-module of rank $\mu$.
We define inductively the sequence $(\underline{a}(k),i(k))\in \NN^{n+1}\times \{0, \ldots ,n\}$ by 
\begin{align*} 
\underline{a}(0)&=(0, \ldots ,0),  & i(0)&=0, \\ 
\underline{a}(k+1)&=\underline{a}(k)+1_{i(k)},  & i(k+1)&=\min\{j\mid \underline{a}(k+1)_{j}/w_{j}\}.
\end{align*} 
For $k\in \{0, \ldots ,\mu-1\}$, we put $
\omega_{k}:=w^{\underline{a}(k_{\min}(s(k)))-\underline{a}(k)}u^{\underline{a}(k)}\omega_{0}
$ where $u^{\underline{a}(k)}:=u_{1}^{\underline{a}(k)_{1}}\ldots
u_{n}^{\underline{a}(k)_{n}}$ and
$w^{\underline{a}(k)}:=w_{1}^{\underline{a}(k)_{1}}\ldots
w_{n}^{\underline{a}(k)_{n}}$.  The classes of $\omega_{0}, \ldots
,\omega_{\mu-1}$ form a $\CC[\theta]$-basis of $G_{0}$, denoted by
$\bs{\omega}$. This basis induces a basis, denoted by $[\bs{\omega}]$,
of the vector space $G_{0}/\theta G_{0}$.  The product structure on
the Jacobian quotient $\mathcal{O}(U)/(\partial f)$ is carried to
$G_{0}/\theta G_{0}$ through the isomorphism $\varphi\mapsto\varphi\omega_{0}$.

\begin{prop}\label{prop:product}
  In the basis $[\bs{\omega}]$ of $G_{0}/\theta G_{0}$, the product is 
  given by 
  \begin{align*} 
  [\omega_{i}]\star[\omega_{j}]=w^{\underline{a}(k_{\min}(s(i)))+\underline{a}(k_{\min}(s(j)))-\underline{a}(k_{\min}(s(\overline{i+j})))+\underline{a}(\overline{i+j})-\underline{a}(i+j)}[\omega_{\overline{i+j}}] 
  \end{align*} 
  where $ \overline{i+j}$ is the sum modulo $\mu$. 
\end{prop} 

\begin{rem}
  This proposition will be used in Section
  \ref{subsec:proof-class-corr} in order to prove the classical
  correspondence. In fact, we will define a product on the graded ring
  $\gr^{\mathcal{N}}_{\star}(G_{0}/\theta G_{0})$ where
  $\mathcal{N}_{\bullet}$ is the Newton filtration of $f$.
\end{rem}

According to Section of $4$ of \cite{DSgm2}, the metric $g$ on
$G_{0}/\theta G_{0}$ in the basis $[\bs{{\omega}}]$ is given by 
Formula (\ref{eq:metric}).

We define 
 \begin{align}\label{eq:3,tenseur}
(\!([a],[b],[c])\!):=g([a]\star [b],[c])
\end{align}
 for any  
$[a],[b],[c]$ in $G_{0}/\theta G_{0}$.
 Proposition \ref{prop:product} and Formula (\ref{eq:metric}) imply 

\begin{cor}\label{cor:tenseur}
  Let  $j,k$ be in $\{0, \ldots ,\mu-1\} $. 
\begin{enumerate} \item If $\overline{1+j+k} \neq n$  then 
  $(\!([{\omega}_{1}],[{\omega}_{j}],[{\omega}_{k}])\!)=0$. 
\item If $\overline{1+j+k} = n$ then   
\begin{align*}(\!([{\omega}_{1}],[{\omega}_{j}],[{\omega}_{k}])\!)= 
 \begin{cases} 
   \prod_{i\in I(j,k)}w_{i}^{-1} & \mbox{ if } 
   \sigma(1)+\sigma(j)+\sigma(k)\neq n\\ 
   \prod_{i\in I(s(j))}w_{i}^{-1} & \mbox{ if } 
   \sigma(1)+\sigma(j)+\sigma(k)=n. 
 \end{cases} 
\end{align*} 
where $I(j,k):=I(s(j))\bigsqcup  I(s(k))$. 
\end{enumerate} 
\end{cor} 

\begin{rem}\label{rem:equivalence}
  \begin{enumerate}
  \item This corollary will be useful to prove the quantum
    correspondence in Section \ref{subsec:proof-quant-corr}.  Because
    the bilinear form $g$ is non degenerated, one can reconstruct the
    multiplication by $[\omega_{1}]$ from this corollary and the
    bilinear form $g$.
  \item Let us remark that the numbers $A_{1jk}(\bs{0})$ in Theorem
    \ref{thm:cond,init} are exactly
    $(\!([{\omega}_{1}],[{\omega}_{j}],[{\omega}_{k}])\!)$.
\end{enumerate}
\end{rem}
 
\subsection{Potential of the Frobenius structure}
In this section we study the potential of the Frobenius structure and
we show that the potential is determined by some numbers (see Theorem \ref{thm:cond,init}). 

 Let $X$ be the base space of a universal unfolding of $f$. Let $t_{0}, \ldots
,t_{\mu-1}$ be the flat coordinates in a neighborhood of $0$ in $X$.
We define the Euler field  by 
\begin{align}\label{eq:champ,euler,sing} 
\mathfrak{E}&=\sum_{k=0}^{\mu-1}(1-\sigma(k))t_{k}\partial 
_{t_{k}}+\mu\partial_{t_{1}}. 
\end{align} 
We develop the potential of the Frobenius structure in series and
we denote it by  
\begin{align*} 
F^{sing}(\bs{t})&=\sum_{\alpha_{0}, \ldots ,\alpha_{\mu-1}\geq 0} 
A(\bs{\alpha}) \frac{t^{\bs{\alpha}}}{\bs{\alpha}!}. 
\end{align*} 
where $\bs{\alpha}:=(\alpha_{0}, \ldots ,\alpha_{\mu-1})$ and
$\frac{t^{\bs{\alpha}}}{\bs{\alpha}!}:=\frac{t_{0}^{\alpha_{0}}}{\alpha_{0}!}\cdots\frac{t_{\mu-1}^{\alpha_{\mu-1}}}{\alpha_{\mu-1}!}$.
We denote $|\alpha|:=\alpha_{0}+\cdots+\alpha_{\mu-1}$ and we call it the
\emph{length} of $\bs{\alpha}$.  We denote by $(g^{ab})$ the inverse matrix
of the non degenerate pairing $g$ in the coordinates $\bs{t}$.  For
any $i,j,k,\ell \in \{0, \ldots ,\mu-1\}$, the potential satisfies the
following conditions :
\begin{align}\label{eq:WDVV} 
(i,j,k,\ell):\ \ \  \sum_{a=0}^{\mu-1}F^{sing}_{ija}
g^{aa^{\star}}F^{sing}_{a^{\star}k\ell}&=\sum_{a=0}^{\mu-1}F^{sing}_{jka}
g^{aa^{\star}}F^{sing}_{a^{\star}i\ell} \mbox{\ \ \  (WDVV equations) } \\
 \mathfrak{E}\cdot
F^{sing}&= (3-n)F^{sing} \mbox{\ \ \ up to  quadratic terms}\label{eq:homo,potentiel} \\
F_{ijk}^{sing}(\bs{0})&=g\mid_{\bs{t}=0}(\partial_{t_{i}}\star \partial_{t_{j}},\partial_{t_{k}}) \mbox{  
  where } F^{sing}_{ijk}:=\frac{\partial^{3} F^{sing}}{\partial 
  t_{i}\partial t_{j} \partial t_{k}}. \label{eq:condition,initiale,potentiel} 
\end{align} 

Note that $g\mid_{\bs{t}=0}(\partial_{t_{i}}\star
\partial_{t_{j}},\partial_{t_{k}})=(\!(\widetilde{\omega}_{i},\widetilde{\omega}_{j},\widetilde{\omega}_{k})\!)$.
Denote $A_{ijk}(\bs{\alpha})$ the number $A(\alpha_{0}, \ldots
,\alpha_{i}+1, \ldots ,\alpha_{j}+1, \ldots ,\alpha_{k}+1, \ldots
,\alpha_{\mu-1})$.

\begin{thm}\label{thm:cond,init}
  The potential $F^{sing}$ is determined by the numbers 
  $A_{1jk}(\bs{0})$ with $j,k\in\{0, \ldots ,\mu-1\}$ such that 
  $\overline{1+j+k}=n$. 
\end{thm} 
 
\begin{rem}\label{rem:basique,potentiel}
  \begin{enumerate}
  \item 
If  $\overline{1+j+k}\neq n$,
  then Condition \eqref{eq:condition,initiale,potentiel} implies that
  $A_{1jk}(\bs{0})=0$. This condition is exactly Proposition
  \ref{prop:invariant,nul} on the A side.
\item If we interpret this theorem on the A side, this means that we have an
algorithm to reconstruct the full quantum cohomology from the small
one.
Note that recently, M.\k Rose has proved in \cite{Rrts} a general
reconstruction theorem for smooth Deligne-Mumford stack. As the small
quantum cohomology of weighted projective spaces are generated by
$H^{2}(|\PP(w)|,\QQ)$ (\cf Corollary 1.2 of \cite{CCLTsqcwps}), Theorem 0.3 of
\cite{Rrts} implies that all genus zero Gromov-Witten invariants can
be reconstructed from the $3$-point invariants. 
\end{enumerate}
\end{rem}

 \begin{proof} 
First, we will show that the potential is determined by the numbers
$A_{ijk}(\bs{0})$ for any $i,j,k\in \{0, \ldots ,\mu-1\}.$

We will show this by induction on the length of the numbers
 $A(\bs{\alpha})$.  For any  $i,j,k,\ell\in \{0, \ldots
   ,\mu-1\}$, the term of
   $F^{sing}_{ija}g^{aa^{\star}}F^{sing}_{a^{\star}k\ell}$ between 
   $\bs{\frac{\bs{t}^{{\alpha}}}{\alpha!}}$ is
   \begin{align*} 
   g^{aa^{\star}}\sum_{\bs{\beta}+\bs{\gamma}=\bs{\alpha}} {\beta_{0} \choose 
     \alpha_{0}}\cdots {\beta_{\mu-1} \choose \alpha_{\mu-1}} 
   A_{ija}(\bs{\beta})A_{a^{\star}k\ell}(\bs{\gamma}). 
   \end{align*} 
   Hence, the terms of the greatest length, that is of length
   $|\bs{\alpha}|+3$, in the sum above are $
   g^{aa^{\star}}A_{ija}(\bs{\alpha})A_{a^{\star}k\ell}(\bs{0})$ and
   $g^{aa^{\star}}A_{ija}(\bs{0})A_{a^{\star}k\ell}(\bs{\alpha})$.  As
   the potential satisfies Conditions
   (\ref{eq:condition,initiale,potentiel}), we deduce that
   $A_{ija}(\bs{0})\neq 0$ if and only if $a=\overline{i+j}^{\star}$.

   In the  WDVV  equation $(1,j,k,\ell)$, the terms of length  
   $|\bs{\alpha}|+3$ in front of  $\frac{t^{\bs{\alpha}}}{\bs{\alpha}!}$ 
   are 
   \begin{align*}
    & g^{\overline{1+j},
       \overline{1+j}^{\star}}A_{1j\overline{1+j}^{\star}}(\bs{0})A_{\overline{1+j}k\ell}(\bs{\alpha}),   &
g^{\overline{k+\ell},
  \overline{k+\ell}^{\star}}A_{1j\overline{k+\ell}}(\bs{\alpha})A_{\overline{k+\ell}^{\star}k\ell}(\bs{0}) ,\\
&g^{\overline{j+k}, \overline{j+k}^{\star}}A_{jk\overline{j+k}^{\star}}(\bs{0})A_{\overline{j+k}1\ell}(\bs{\alpha})    ,&
g^{\overline{1+\ell}, \overline{1+\ell}^{\star}}A_{jk\overline{1+\ell}}(\bs{\alpha})A_{\overline{1+\ell}^{\star}1\ell}(\bs{0}). 
\end{align*}
The terms of the form  $A_{???}(\bs{0})$ are computed by  
 (\ref{eq:condition,initiale,potentiel}) and the homogeneity condition
  (\ref{eq:homo,potentiel}) implies that  
 \begin{align*} 
 A(\alpha_{0},\alpha_{1}+1,\alpha_{2},\ldots,\alpha_{\mu-1})=\frac{1}{\mu}A(\bs{\alpha})d(\bs{\alpha}) 
 \mbox{ pour } |\bs{\alpha}|\geq 3 
 \end{align*} 
 where $d(\bs{\alpha})=3-n+\sum_{k=0}^{\mu-1}\alpha_{k}(\sigma(k)-1)$. 
 Hence, we can express the numbers $A_{1??}(\bs{\alpha})$ with numbers
 of length strictly smaller.  The 
 WDVV equation $(1,j,k,\ell)$ gives a relation between  
 $A_{\overline{1+j}k\ell}(\bs{\alpha})$ and 
 $A_{jk\overline{1+\ell}}(\bs{\alpha})$. We conclude that after a
 finite number of steps, we can express any number
 $A_{ijk}(\bs{\alpha})$ with terms of strictly smaller length. By
 induction we have that the potential is determined by the numbers
$A_{ijk}(\bs{0})$ for any $i,j,k\in \{0, \ldots ,\mu-1\}.$

In order to finish the proof, it is enough to show that
 we can compute any numbers $A_{ijk}(\bs{0})$ from the numbers of the
 form $A_{1??}(\bs{0})$.  The numbers of length $3$ in  
   the equation $(1,j,k,\ell)$ are non zero if and only if  
   $\overline{1+j+k+\ell}=n$.  Under this condition, we have  
   $\overline{1+j}=\overline{k+\ell}^{\star}$ and
   $\overline{j+k}^{\star}=\overline{1+\ell}$.  Hence, the terms of length 
    $3$ in the equation $(1,j,k,\ell)$ are 
    \begin{align*}
       A_{1j\overline{1+j}^{\star}}(\bs{0})A_{\overline{1+j}k\ell}(\bs{0})
 \hspace{2cm}\mbox{ and }&&
   A_{jk\overline{1+\ell}}(\bs{0})A_{\overline{1+\ell}^{\star}1\ell}(\bs{0}).
    \end{align*}
Considering successively the WDVV equations  
 $(1,j-1,k,\ell+1)$,$(1,j-2,k,\ell+2)$,..., we can express  
 $A_{\overline{1+j}k\ell}(\bs{0})$ in terms of the numbers of the form 
 $A_{1??}(\bs{0})$. 
\end{proof}

\section{Correspondences}

\subsection{Combinatorics of numbers $\sigma$}
\label{subsec:Comb-numb-sigma}

We define an order on the circle $S^{1}$ by choosing the principal determination of
the argument. Choose a non decreasing bijection $\widetilde{s}:\{0, \ldots
,\mu-1\}\to \sqcup \bs{\mu}_{w_{i}}$.
For any $g\in \sqcup \bs{\mu}_{w_{i}}$, we put
\begin{align*}
  k_{\min}(g):=\min\{i\in\{0, \ldots ,\mu-1\} \mid \widetilde{s}(i)=g\}.
\end{align*}
For any $g\in \sqcup\bs{\mu}_{w_{i}}$, we denote $\gamma(g)$ the unique
element in $[0,1[$ such that $\exp(2i\pi\gamma(g))=g$.
Recall that we have defined $\sigma(\cdot)$ in Notation \ref{not:sigma}. We have that
$\sigma(i)=i-\mu\gamma(\widetilde{s}(i))$.

The following proposition is straightforward.
\begin{prop}\label{prop:comb}
  \begin{enumerate}
  \item\label{item:prop:comb,1} For any $g\in  \sqcup\bs{\mu}_{w_{i}}$, we have
 \begin{align*}
  k_{\min}(g)= \codim\PP(w)_{(g)}+\sum [\gamma(g)w_{i}]
\end{align*}where $[\cdot]$ means the integer part.
\item\label{item:prop:comb,2} For any $g\in \sqcup \bs{\mu}_{w_{i}}$ and $d\in\{0, \ldots
  ,\dim\PP(w)_{(g)}\}$, we have 
  \begin{align*}
   2\sigma(k_{\min}(g^{-1})+d)=\deg(\eta^{d}_{g}). 
  \end{align*}
\item\label{item:prop:comb,3} For any $g,g'\in \sqcup \bs{\mu}_{w_{i}}$ and $(d,d')\in\{0, \ldots
  ,\dim\PP(w)_{(g)}\}\times \{0, \ldots
  ,\dim\PP(w)_{(g')}\}$, we have the following equivalence
  \begin{align*}
k_{\min}(g^{-1})+d+k_{\min}(g'^{-1})+d'=n \mod \mu \\
\Leftrightarrow
\begin{cases}
  g.g'=\id\\
d+d'=\dim\PP(w)_{(g)}.
\end{cases}
  \end{align*}
  \end{enumerate}
\end{prop}

\subsection{Proof of the classical correspondence}\label{subsec:proof-class-corr}
Denoted by $\mathcal{N}_{\bullet} G$, the Newton filtration of  the
Gauss-Manin system $G$ (see Paragraph $2.$e of \cite{DSgm1}). This
filtration induces a filtration on $G_{0}/\theta G_{0}$, denoted by
$\mathcal{N}_{\bullet}(G_{0}/\theta G_{0})$ or just $\mathcal{N}_{\bullet}$ when there is
no ambiguity.
Let $\Xi$ be the $\CC$-linear map defined by 
 \begin{align*}
    \Xi \,: H^{2\star}_{\orb}(\PP(w),\CC) & \longrightarrow 
    \gr^{\mathcal{N}}_{\star}\left(G_{0}/\theta G_{0}\right) \\ 
    \eta^{d}_{g} & \longmapsto 
    [\![ {\omega}_{k_{\min}(g^{-1})+d}]\!] \nonumber 
    \end{align*} 
    where $[\![\cdot]\!]$ means the class in
    $\gr^{\mathcal{N}}_{\star}\left(G_{0}/\theta G_{0}\right)$.  The
    non-degenerate bilinear form $g(\cdot,\cdot)$ on $G_{0}/\theta
    G_{0}$ induces a non-degenerate bilinear form, denoted by
    $[\![g]\!](\cdot,\cdot)$ on
    $\gr^{\mathcal{N}}_{\star}\left(G_{0}/\theta G_{0}\right)$.
    Because for any $\beta_{1},\beta_{2}\in \QQ$, we have that
    $\mathcal{N}_{\beta_{1}}(G_{0}/\theta G_{0}) \star
    \mathcal{N}_{\beta_{2}}(G_{0}/\theta G_{0})$ is included in
    $\mathcal{N}_{\beta_{1}+\beta_{2}}(G_{0}/\theta G_{0})$ (see
    Proposition VI.$3.1$ of \cite{Mthesis}), we can define a product,
    denoted by $\cup$, on $\gr^{\mathcal{N}}_{\star}\left(G_{0}/\theta
      G_{0}\right)$.

 \begin{thm}\label{thm:corr,classique} The map $\Xi$ is a graded isomorphism   
   between the graded Frobenius algebras $
   (H^{2\star}_{\orb}(\PP(w),\CC),\cup,\langle\cdot,\cdot\rangle) $
   and $ (\gr^{\mathcal{N}}_{\star} \left(G_{0}/\theta
     G_{0}\right),\cup,[\![g]\!](\cdot,\cdot)). $
 \end{thm}

\begin{proof}  
  According to Corollary \ref{cor:degre,coho,orbifold}, we have $
  \deg(\eta^{d}_{\gamma})=2(d+\age(g))$. 
Proposition
 \ref{prop:comb}.\eqref{item:prop:comb,2} implies that $\Xi(\eta^{d}_{\gamma})$ is in the graded
  $\gr^{\mathcal{N}}_{d+a(\gamma)}\left(G_{0}/\theta G_{0}\right)$.
  We conclude that $\Xi$ is a graded map.  On one hand Proposition
  \ref{prop:dualite,poincare}, Formula (\ref{eq:metric}) and
  Proposition \ref{prop:comb}.\eqref{item:prop:comb,3}  imply
  that $\langle
  \eta^{d}_{\gamma},\eta^{d'}_{\gamma'}\rangle=[\![g]\!](\Xi(\eta^{d}_{\gamma}),\Xi(\eta^{d'}_{\gamma'})).$
  On the other hand Corollary \ref{cor:cup,produit} and Proposition
  \ref{prop:product} imply that $
  \Xi(\eta^{d_{0}}_{\gamma_{0}}\cup\eta^{d_{1}}_{\gamma_{1}})=\Xi(\eta^{d_{0}}_{\gamma_{0}})\cup\Xi(\eta^{d_{1}}_{\gamma_{1}}).$
\end{proof}

\subsection{Proof of the quantum correspondence}\label{subsec:proof-quant-corr} 
 
 Let $\widetilde{\Xi}$ be the $\CC$-linear map defined by 
 \begin{align*}
   \widetilde{\Xi} \,: QH^{2\star}_{\orb}(\PP(w),\CC) &
   \longrightarrow
   G_{0}/\theta G_{0} \\
   \eta_{g}^{d} & \longmapsto [ {\omega}_{k_{\min}(g^{-1})+d}] \end{align*}
 
This map is an isomorphism of vector space.
\begin{cor} Let $w_{0}, \ldots 
  ,w_{n}$ be integers. The Frobenius manifolds associated to
  the Laurent polynomial $f$ and the orbifold $\PP(w)$ are isomorphic. 
 \end{cor}  

\begin{proof} 
As the Frobenius manifold associated to $f$ is semi simple (\cf Theorem
\ref{thm:cond,init,sing}) it is enough to show that both Frobenius manifolds carry the same
  initial conditions.
 
  Theorem \ref{thm:cond,init,sing} gives us the initial condition
  $(A_{0}^{\circ},A_{\infty},e_{0},g)$ for the Frobenius manifold
  associated to the Laurent polynomial $f$.  Theorem
  \ref{thm:corr,classique} implies that we
  have the same matrix $A_{\infty}$, the same eigenvector $e_{0}$ for
  the eigenvalue $q=0$ and the same bilinear non-degenerate form.
  
  We have to compare the matrices $A_{0}^{\circ}$ which correspond to
  the multiplication by the Euler fields at the origin.  Formulas
  (\ref{eq:euler,field,Aside}) and \eqref{eq:champ,euler,sing} show
  that the Euler fields are the same.  Corollary \ref{cor:3,tenseur},
  via Propositions \ref{prop:invariant,nul} and
  \ref{prop:invariants,durs}, and Corollary \ref{cor:tenseur} imply
  that for any $g,g'\in \cup \bs{\mu}_{w_{i}}$ and $(d,d')\in \{0,
  \ldots ,\dim\PP(w)_{g}\}\times\{0, \ldots ,\PP(w)_{(g')}\} $, we
  have
    \begin{align*} 
      \langle\eta_{1}^{1}\star\eta_{g}^{d},\eta_{g'}^{d'}
\rangle&=(\!(\eta_{1}^{1},\eta_{g}^{d},\eta_{g'}^{d'})\!)\\
&=(\!(\widetilde{\Xi}(\eta_{1}^{1}),\widetilde{\Xi}(\eta_{g}^{d}),\widetilde{\Xi}(\eta_{g'}^{d'}))\!)\\
&=  g(\widetilde{\Xi}(\eta_{1}^{1})\star\widetilde{\Xi}(\eta_{g}^{d}),\widetilde{\Xi}(\eta_{g'}^{d'}))
\end{align*}
Moreover, Proposition
  \ref{prop:dualite,poincare} and Formula (\ref{eq:metric}) imply
  that we have $\langle \eta_{g}^{d},\eta_{g'}^{d'}\rangle=g(\widetilde{\Xi}(\eta_{g}^{d}),\widetilde{\Xi}(\eta_{g'}^{d'})).$
Hence, we deduce  that the multiplications by the Euler field at the origin are the
same. 
\end{proof}  

\section{Appendix : small quantum cohomology of weighted projective spaces} 
\label{appen}

In \cite{CCLTsqcwps}, T.\k Coates, A.\k Corti, Y.-P. \k Lee and
H.-H.\k Tseng has computed the small quantum cohomology of weighted
projective spaces. We recall their results with our notation.

Denote the elements of $\cup
\mu_{w_{i}}$ by $1=g_{0}<\cdots<g_{\delta}$  where the order is
defined by choosing the principal determination of the argument (\cf
Notation \ref{not:order}). Recall that for any $g\in \cup
\bs{\mu}_{w_{i}}$, there exists a unique $\gamma(g)\in[0,1[$ such that
$g=\exp(2i\pi\gamma(g))$. 

For any $k\in\{0, \ldots ,\delta\}$, put 
\begin{displaymath}
  s_{k}=
  \begin{cases}
    1 & \mbox{ if } k=1 \\
\frac{\ds{\prod_{\gamma(g_{m})<\gamma(g_{k})}(\gamma(g_{k})-\gamma(g_{m}))^{\dim\PP(w)_{(g_{m})}+1}}}{\ds{\prod_{i=0}^{n}(\gamma(g_{j})w_{i})^{\underline{\lceil
      \gamma(g_{k})w_{i}    \rceil}}}}
& \mbox{ otherwise}
  \end{cases}
\end{displaymath}
where we put 
\begin{align*}
  x^{\underline{n}}=x(x-1)\ldots(x-n+1).
\end{align*}

According to Corollary 1.2 of \cite{CCLTsqcwps}, the small quantum cohomology of weighted
projective space is generated by $\eta_{1}^{1}$ and $\eta_{g}^{0}$
for any $g\in \cup \bs{\mu}_{w_{i}}$. The relations are
\begin{align}
  (\eta_{1}^{1})^{k_{\min}(g_{k})}&=Q^{\gamma(g_{k})}s_{k}\eta_{g_{k}^{-1}}^{0} \label{eq:5}\\
 \underbrace{\eta_{1}^{1}\star \cdots \star\eta_{1}^{1}}_{\sum w_{i}}&= Q\prod w_{i}^{-w_{i}}\nonumber
\end{align}
where $k_{\min}(g):=\sum [w_{i}\gamma(g)]+\codim \PP(w)_{(g)}$ (\cf
Section $6.a$ for an other interpretation of $k_{\min}$) and $Q$ is a
formal variable of degre $\mu$.  The careful reader will notice that
Equation (\ref{eq:5}) is not exactly the one of Corollary 1.2 in
\cite{CCLTsqcwps}. Indeed, they define weighted projective spaces
differently that is as the
quotient stack of $[\CC^{n+1}-\{0\}/\CC^{\star}]$ where
$\CC^{\star}$ acts with weights $-w_{0}, \ldots ,-w_{n}$.

\begin{prop}\label{prop:mult,eta,1,1}
  For any $k\in\{0, \ldots ,\delta\}$, we have that
  \begin{displaymath}
    \eta_{1}^{1}\star\eta_{g_{k-1}^{-1}}^{\dim\PP(w)_{(g_{k-1})}}
= \eta_{g_{k}^{-1}}^{0}Q^{\gamma(g_{k})-\gamma(g_{k-1})}\prod_{i\in
      I(g_{k-1})}w_{i}^{-1}.
  \end{displaymath}
\end{prop}

\begin{proof}
  First, we will show that for any $k\in\{0, \ldots ,\delta\}$, we
  have 
  \begin{align}\label{eq:6}
    s_{k}=\prod_{i=0}^{n} w_{i}^{-\lceil\gamma(g_{k})w_{i}\rceil}.
  \end{align}
  For any $k\in\{0, \ldots ,\delta\}$ and any $i\in\{0, \ldots ,n\}$, we have 
\begin{align*}
\frac{\lceil\gamma(g_{k})w_{i}\rceil-1}{w_{i}}<
\gamma(g_{k})\leq\frac{\lceil\gamma(g_{k})w_{i}\rceil}{w_{i}}.
\end{align*}
We deduce that 
\begin{align*}
  \prod_{i=0}^{n}(\gamma(g_{j})w_{i})^{\underline{\lceil
      \gamma(g_{k})w_{i}    \rceil}} & = w_{i}^{\lceil
    \gamma(g_{k})_{i}\rceil}\hspace{-.5cm}\prod_{\stackrel{\ell \ \mid}{\ell/w_{i}<\gamma(g_{k})}}
  \hspace{-.5cm}\left( \gamma(g_{k})-\frac{\ell}{w_{i}}\right).
\end{align*}
We deduce Formula (\ref{eq:6}).

Put  $d(g_{k}):=\dim\PP(w)_{(g_{k})}$.
 We have
 \begin{align*}
  \eta_{g_{k-1}^{-1}}^{d(g_{k-1})} = \eta_{g_{k-1}^{-1}}^{0}\star (\eta_{1}^{1})^{d(g_{k-1})}.
\end{align*}
    Section $6.a$ and Proposition $6.1$ imply  that
    $k_{\min}(g_{k-1})+d(g_{k-1})+1=k_{\min}(g_{k})$. 
 We deduce that
\begin{align*}
   \eta_{1}^{1} \star \eta_{g_{k-1}^{-1}}^{d(g_{k-1})} &=
   (\eta_{1}^{1})^{k_{\min}(g_{k})}Q^{-\gamma(g_{k-1})} \prod_{i=0}^{n} w_{i}^{\lceil
    \gamma(g_{k-1})w_{i}\rceil} \\
&=\eta_{g_{k}^{-1}}^{0} Q^{\gamma(g_{k})-\gamma(g_{k-1})}\prod_{i=0}^{n} w_{i}^{\lceil
     \gamma(g_{k-1})w_{i}\rceil-\lceil\gamma(g_{k})w_{i}\rceil} \\
\end{align*}
Then the following lemma finishes the proof.
\end{proof}

\begin{lem}
  For any $i\in \{0, \ldots ,n\}$, we have 
  \begin{align*}
    \lceil w_{i}\gamma(g_{k-1})\rceil -\lceil w_{i}\gamma(g_{k})\rceil = 
    \begin{cases}
      -1 & \mbox{ if } i\in I(g_{k-1})\\
  0 & \mbox{ otherwise.}
    \end{cases}
  \end{align*}
\end{lem}

\begin{proof}
  By choosing the principal determination of the argument, we order
  the elements in $\cup \bs{\mu}_{w_{i}}$ by
  $1=g_{0}<g_1<\cdots<g_{d}$.  Hence we have
\begin{align*}
  0<w_{i}\gamma(g_{1})<\cdots<1<\cdots<2<\cdots<w_{i}-1<\cdots<w_{i}\gamma(g_{d}).
\end{align*}
The above formula implies the following alternative :
\begin{itemize}
\item if $w_{i}\gamma(g_{k-1})\in \NN$ (\ie $i\in I(g_{k-1})$), we
  have $\lceil w_{i}\gamma(g_{k-1})\rceil -\lceil
  w_{i}\gamma(g_{k})\rceil =-1$.
\item if $w_{i}\gamma(g_{k})\in \NN$, we have $\lceil
  w_{i}\gamma(g_{k-1})\rceil -\lceil w_{i}\gamma(g_{k})\rceil =0$.
\item if $w_{i}\gamma(g_{k-1}),w_{i}\gamma(g_{k})\notin \NN$, we
  have $\lceil w_{i}\gamma(g_{k-1})\rceil -\lceil w_{i}\gamma(g_{k})\rceil
  =0$.
\end{itemize}
\end{proof}

\bibliographystyle{alpha} 
\bibliography{biblio} 

\begin{thebibliography}{CCLT06}

\bibitem[AGV06]{AGVgwdms}
Dan Abramovich, Tom Graber, and Angelo Vistoli.
\newblock {Gromov--Witten theory of Deligne--Mumford stacks}.
\newblock {\em math.AG/0603151}, page~57, 2006.

\bibitem[Bar00]{Bms}
Serguei Barannikov.
\newblock {Semi-infinite Hodge structures and mirror symmetry for projective
  spaces}.
\newblock {\em Math.AG/0010157}, page~17, 2000.

\bibitem[BCS05]{BCSocdms}
Lev~A. Borisov, Linda Chen, and Gregory~G. Smith.
\newblock The orbifold {C}how ring of toric {D}eligne-{M}umford stacks.
\newblock {\em J. Amer. Math. Soc.}, (18)(1):193--215 (electronic), 2005.

\bibitem[CCLT06]{CCLTsqcwps}
Tom Coates, Alessio Corti, Y.-P. Lee, and Hsian-Hua Tseng.
\newblock {Small quantum orbifold cohomology of weighted projective spaces}.
\newblock {\em math.AG/0608481}, page~50, 2006.

\bibitem[CH04]{CH}
Bohui Chen and Shengda Hu.
\newblock {A deRham model for Chen-Ruan cohomology ring of abelian orbifolds}.
\newblock {\em math.SG/0408265}, page~13, 2004.

\bibitem[CR02]{CRogw}
Weimin Chen and Yongbin Ruan.
\newblock Orbifold {G}romov-{W}itten theory.
\newblock In {\em Orbifolds in mathematics and physics (Madison, WI, 2001)},
  volume (310) of {\em Contemp. Math.}, pages 25--85. Amer. Math. Soc.,
  Providence, RI, 2002.

\bibitem[CR04]{CRnco}
Weimin Chen and Yongbin Ruan.
\newblock A new cohomology theory of orbifold.
\newblock {\em Comm. Math. Phys.}, (248)(1):1--31, 2004.

\bibitem[DS03]{DSgm1}
Antoine Douai and Claude Sabbah.
\newblock Gauss-{M}anin systems, {B}rieskorn lattices and {F}robenius
  structures. {I}.
\newblock {\em Ann. Inst. Fourier (Grenoble)}, (53)(4):1055--1116, 2003.

\bibitem[DS04]{DSgm2}
Antoine Douai and Claude Sabbah.
\newblock Gauss-{M}anin systems, {B}rieskorn lattices and {F}robenius
  structures. {II}.
\newblock In {\em Frobenius manifolds}, Aspects Math., E36, pages 1--18.
  Vieweg, Wiesbaden, 2004.

\bibitem[Dub96]{Dtft}
Boris Dubrovin.
\newblock Geometry of {$2$}{D} topological field theories.
\newblock In {\em Integrable systems and quantum groups (Montecatini Terme,
  1993)}, volume (1620) of {\em Lecture Notes in Math.}, pages 120--348.
  Springer, Berlin, 1996.

\bibitem[FO99]{FOgwi}
Kenji Fukaya and Kaoru Ono.
\newblock Arnold conjecture and {G}romov-{W}itten invariant for general
  symplectic manifolds.
\newblock In {\em The Arnoldfest (Toronto, ON, 1997)}, volume (24) of {\em
  Fields Inst. Commun.}, pages 173--190. Amer. Math. Soc., Providence, RI,
  1999.

\bibitem[GH94]{GHag}
Phillip Griffiths and Joseph Harris.
\newblock {\em Principles of algebraic geometry}.
\newblock Wiley Classics Library. John Wiley \& Sons Inc., New York, 1994.
\newblock Reprint of the 1978 original.

\bibitem[Her02]{Hfm}
Claus Hertling.
\newblock {\em Frobenius manifolds and moduli spaces for singularities}, volume
  (151) of {\em Cambridge Tracts in Mathematics}.
\newblock Cambridge University Press, Cambridge, 2002.

\bibitem[HM04]{HMumc}
Claus Hertling and Yuri Manin.
\newblock Unfoldings of meromorphic connections and a construction of
  {F}robenius manifolds.
\newblock In {\em Frobenius manifolds}, Aspects Math., E36, pages 113--144.
  Vieweg, Wiesbaden, 2004.

\bibitem[Jia03]{Jocwps}
Yunfeng Jiang.
\newblock {The Chen-Ruan cohomology of weighted projective spaces}.
\newblock {\em math.AG/0304140}, page~34, 2003.

\bibitem[Kaw73]{Kcwps}
Tetsuro Kawasaki.
\newblock Cohomology of twisted projective spaces and lens complexes.
\newblock {\em Math. Ann.}, (206):243--248, 1973.

\bibitem[KM94]{KMgwqceg}
Maxim Kontsevich and Yuri Manin.
\newblock Gromov-{W}itten classes, quantum cohomology, and enumerative
  geometry.
\newblock {\em Comm. Math. Phys.}, (164)(3):525--562, 1994.

\bibitem[LMB00]{LMBca}
G{\'e}rard Laumon and Laurent Moret-Bailly.
\newblock {\em Champs alg\'ebriques}, volume~39 of {\em Ergebnisse der
  Mathematik und ihrer Grenzgebiete. 3. Folge. A Series of Modern Surveys in
  Mathematics [Results in Mathematics and Related Areas. 3rd Series. A Series
  of Modern Surveys in Mathematics]}.
\newblock Springer-Verlag, Berlin, 2000.

\bibitem[Man99]{Mfm}
Yuri Manin.
\newblock {\em Frobenius manifolds, quantum cohomology, and moduli spaces},
  volume (47) of {\em American Mathematical Society Colloquium Publications}.
\newblock American Mathematical Society, Providence, RI, 1999.

\bibitem[Man05]{Mthesis}
Etienne Mann.
\newblock {C}ohomologie quantique orbifolde des espaces projectifs \`a poids.
\newblock {\em math.AG/0510331}, page 136, 2005.

\bibitem[MP97]{MPos}
Ieke Moerdijk and Dorette~A. Pronk.
\newblock Orbifolds, sheaves and groupoids.
\newblock {\em $K$-Theory}, (12)(1):3--21, 1997.

\bibitem[Ros06]{Rrts}
A.Michael Rose.
\newblock {A reconstruction theorem for genus zero Gromov-Witten invariants of
  stacks}.
\newblock {\em math.AG/0605776}, page~14, 2006.

\bibitem[Sab02]{Sdivf}
Claude Sabbah.
\newblock {\em D\'eformations isomonodromiques et vari\'et\'es de {F}robenius}.
\newblock Savoirs Actuels (Les Ulis). [Current Scholarship (Les Ulis)]. EDP
  Sciences, Les Ulis, 2002.
\newblock Math\'ematiques (Les Ulis). [Mathematics (Les Ulis)].

\bibitem[Sat56]{Sgm}
Ichir{\^o} Satake.
\newblock On a generalization of the notion of manifold.
\newblock {\em Proc. Nat. Acad. Sci. U.S.A.}, (42):359--363, 1956.

\bibitem[Sat57]{Sgb}
Ichir{\^o} Satake.
\newblock The {G}auss-{B}onnet theorem for {$V$}-manifolds.
\newblock {\em J. Math. Soc. Japan}, (9):464--492, 1957.

\end{thebibliography}
\end{document}